\documentclass[12pt]{article}
\usepackage{amssymb,amsfonts,amsmath}
\usepackage{epsfig}
\usepackage{xcolor}
\usepackage{graphicx}
\usepackage{float,wrapfig}

\numberwithin{equation}{section}
\numberwithin{table}{section}
\numberwithin{figure}{section}

\newcommand{\be}{\begin{equation}}
\newcommand{\ee}{\end{equation}}
\newcommand{\CO}{\mbox{$\mathrm{CO}_2$}}
\newcommand{\R}{\mathbb{R}}
\newcommand{\Z}{\mathbb{Z}}
\newcommand{\sech}{\mbox{$\mathrm{sech}$}}
                       
\title{Modeling the Dynamics of Glacial Cycles}

\author{
Hans Engler${}^{1, 2}$, 
Hans G. Kaper${}^{1, 2}$,
Tasso J. Kaper${}^{3}$, 
Theodore Vo${}^{3}$ \\
${}^1$ Department of Mathematics and Statistics, \\
 Georgetown University, Washington, DC \\
${}^2$ Mathematics and Climate Research Network (MCRN) \\
${}^3$ Department of Mathematics and Statistics, \\
 Boston University, Boston, MA
}

\date{ }

\begin{document}
\bibliographystyle{siam}

\maketitle

\abstract{
This article is concerned with the dynamics of glacial cycles
observed in the geological record of the Pleistocene Epoch.
It focuses on a conceptual model proposed by Maasch and Saltzman
[J.\ Geophys.\ Res., \textbf{95}, D2 (1990), pp.\ 1955--1963],
which is based on physical arguments and emphasizes
the role of atmospheric \CO\ in the generation and persistence
of periodic orbits (limit cycles).
The model consists of
three ordinary differential equations with four parameters for
the anomalies of the total global ice mass, 
the atmospheric \CO\ concentration, and
the volume of the North Atlantic Deep Water.
In this article, it is shown that a simplified two-dimensional
symmetric version displays many of the essential features of the full model,
including equilibrium states, limit cycles, their basic bifurcations,
and a Bogdanov--Takens point that serves as an organizing center
for the local and global dynamics.
Also, symmetry breaking splits the Bogdanov--Takens point into two,
with different local dynamics in their neighborhoods.
}

\section{Introduction\label{s-intro}}
Earth's climate during the \emph{Pleistocene Epoch}---the geological period
from approximately 2.6~million years before present (2.6\,Myr~BP)
until approximately 11.7~thousand years before present (11.7\,Kyr~BP)---is of
great interest in the geosciences community.
The geological record gives evidence of cycles of advancing
and retreating continental glaciers and ice sheets,
mostly at high latitudes and high altitudes
and especially in the Northern Hemisphere.

To reconstruct the Pleistocene climate, geoscientists rely on
geological proxies, particularly a dimensionless quantity 
denoted by $\delta^{18} \mathrm{O}$, 
which is measured in parts per mille.
This quantity measures the deviation of the ratio
${}^{18}\mathrm{O} / {}^{16}\mathrm{O}$ 
of the stable oxygen isotopes 
${}^{18}\mathrm{O}$ and ${}^{16}\mathrm{O}$
in a given marine sediment sample
from the same ratio in a universally accepted standard sample. 
The relative amount of the isotope ${}^{18}\mathrm{O}$ in ocean water
is known to be higher at tropical latitudes than near the poles,
since water with the heavier oxygen isotope is slightly less likely
to evaporate and more likely to precipitate first.
Similarly, water with the lighter isotope ${}^{16}\mathrm{O}$
is more likely to be found in ice sheets and in rain water
at high latitudes, 
since it is favored in atmospheric transport across latitudes. 
The global distribution of $\delta^{18}\mathrm{O}$ in ocean water
therefore varies in a known way between glacial and interglacial periods. 
A record of these variations is preserved in the calcium carbonate shells
of foraminifera, a class of common single cell marine organisms. 
These fossil records may be sampled from deep sea sediment cores,
and their age and $\delta^{18} \mathrm{O}$ may be determined. 
Details are described in~{lisiecki2005pliocene}. 

\begin{figure}[ht]
\begin{center}
\resizebox{0.9\textwidth}{!}{\includegraphics{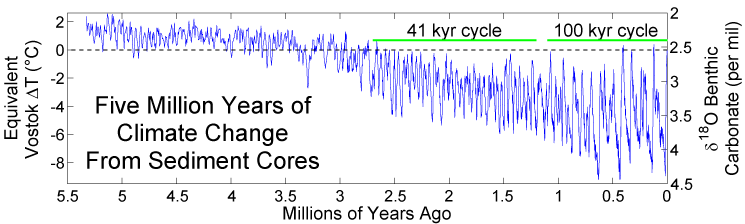}}
\caption{Time series of the isotope ratio 
$\delta^{18} \mathrm{O} = {}^{18} \mathrm{O} / {}^{16} \mathrm{O}$
(scale on the right, in parts per mille) for the past 5.3\,million years~\cite{lisiecki2005pliocene}.
The scale on the left gives the equivalent temperature anomaly (in degrees Celsius).
[Source: Wikipedia, Marine Isotope Stage]
\label{f-TempRecord5Myr}}
\end{center}
\end{figure}

Figure~\ref{f-TempRecord5Myr} shows the LR04 time series of $\delta^{18} \mathrm{O}$
over the past 5.3\,million years, reconstructed from sediment core data
collected at 57 geographically distributed sites around the globe~\cite{lisiecki2005pliocene}.
As the observed isotope variations are similar in shape to the temperature variations
reconstructed from ice core data for the past 420\,Kyr at the Vostok Station in Antarctica,
the values of $\delta^{18} \mathrm{O}$ (right scale) have been aligned with
the reported temperature variations from the Vostok ice core (left scale)~\cite{Petit1999}.
The graph shows a relatively stable temperature during the period preceding the Pleistocene
and increasing variability during the Pleistocene.

\begin{figure}[ht]
\begin{center}
\resizebox{0.9\textwidth}{!}{\includegraphics{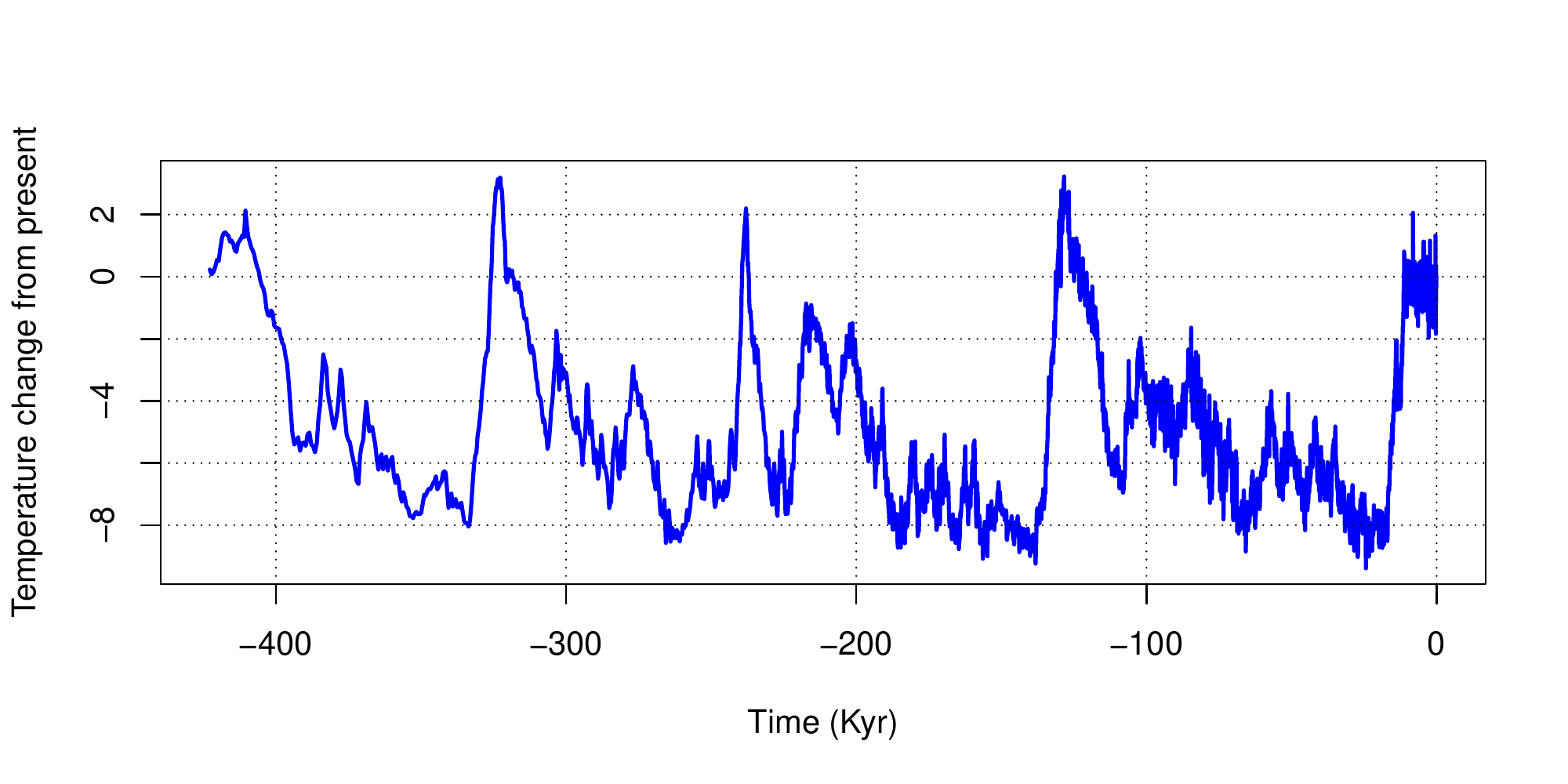}}
\caption{Time series of the global mean temperature for the past 420,000 years.
[Data from Carbon Dioxide Information Analysis Center (CDIAC) 
at Oak Ridge National Laboratory]
\label{f-TempRecord400Kyr}
}
\end{center}
\end{figure}

The typical pattern throughout most of the Pleistocene
resembles that of a sawtooth wave, where a slow glaciation
is followed by a rapid deglaciation.
In the early Pleistocene, 
until approximately 1.2\,Myr~BP, 
the period of a glacial cycle averages 41\,Kyr;
after the mid-Pleistocene transition,
which occurred from approximately 1.2\,Myr~BP 
until approximately 0.8\,Myr~BP,
the glacial cycles of the late Pleistocene
have a noticeably greater amplitude,
and their period averages 100\,Kyr. 
Figure~\ref{f-TempRecord400Kyr} shows
the global mean temperature for the past
420\,Kyr, reconstructed from Vostok ice core data.
The 100\,Kyr cycle and the sawtooth pattern
are clearly visible.

These observations suggest a number of questions for climate science. 
What caused the glacial oscillations during the Pleistocene?
Why were the periods of the glacial cycles during the early
and late Pleistocene different?
What could possibly have caused the transition 
from 41\,Kyr cycles to 100\,Kyr cycles during 
the mid-Pleistocene?

In this article, we discuss a conceptual model of the Pleistocene climate
proposed by Maasch and Saltzman in~\cite{MaaschSaltzman1990}
to explain the phenomenon of glacial cycles.
The model is conceptual, in the sense that it describes the state of the climate
in a few variables, ignoring most of the processes that go into a complete
climate model, but still captures the essence of the phenomenon.
It is based on sound physical principles and, as we will see,
makes for an interesting application of dynamical systems theory.

The numerical continuation results for the bifurcation curves
reported in this article were obtained
using the software package AUTO~\cite{DCFKOPSWZ2007}; 
see also~\cite{D1981, DKK1991}.
Some recent texts on issues of climate dynamics
are~{Cook2013, Dijkstra2013, OT131, MarshallPlumb2007, SunBryan2013}.

In Section~\ref{s-Background}, we present background information
to motivate the particular choices underlying the Maasch--Saltzman model.
In Section~\ref{s-MSModel}, we derive the Maasch--Saltzman model
from physical principles and formulate it as a dynamical system
in a three-dimensional state space with four parameters.
In Section~\ref{s-SimplifyMS}, we introduce two simplifications
that render the Maasch--Saltzman model symmetric and reduce
it to a two-dimensional dynamical system with two parameters
that can be analyzed rigorously and completely.
In Section~\ref{s-2Dasym}, we introduce asymmetry into the simplified
two-dimensional model and show the effects of symmetry breaking.
In the final Section~\ref{s-Summary}, we summarize our results.

\section{Background\label{s-Background}}

There is general agreement that the periodicity of the glacial cycles
is related to variations in the Earth's orbital parameters~\cite{JGRC:JGRC3036}.
To first order, Earth's climate is driven by the Sun.
The Earth receives energy from the Sun in the form
of ultraviolet (short wavelength) radiation.
This energy is redistributed around the globe
and eventually reemitted into space in the form
of infrared (long wavelength) outgoing radiation.
The amount of energy reaching the top of the atmosphere
per unit area and per unit time is known as the \emph{insolation}
(\emph{in}cident \emph{sol}ar radi\emph{ation}),
which is measured in watts per square meter.

\subsection{Orbital Forcing\label{ss-OrbitalForcing}}
The insolation varies with the distance from the Earth to the Sun
and thus depends on Earth's orbit around the Sun.
This is the basis of the Milankovitch theory of 
\emph{orbital forcing}~\cite{OT131, Milankovic1941}.

\begin{minipage}{0.4\textwidth}
\resizebox{\textwidth}{!}{\includegraphics{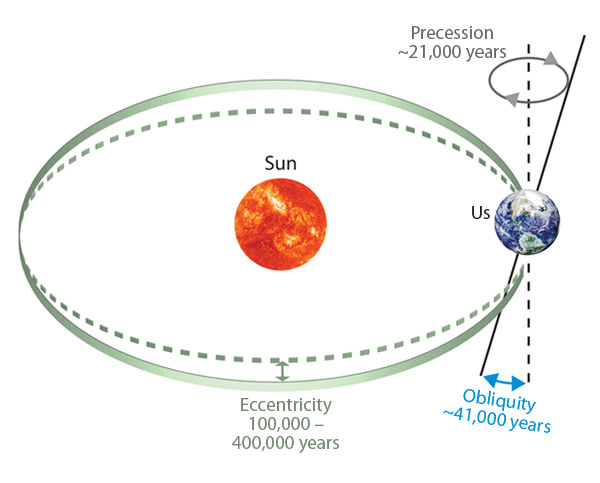}}
\end{minipage}
\hfill
\begin{minipage}{0.4\textwidth}
The Earth moves around the Sun in an elliptical orbit;
its \emph{eccentricity} varies with time but has dominant frequencies
at approximately 100\,Kyr and 400\,Kyr.
As the Earth moves around the Sun, it rotates around its axis.
The axis is tilted with respect to the normal to the orbital plane;
the tilt, known as \emph{obliquity}, is also close to periodic, 
with a dominant frequency of approximately 41\,Kyr.
(This tilt is the main cause of the seasonal variation of our climate.)
\end{minipage}

\vspace{0.6ex} \noindent  
In addition, the Earth is like a spinning top wobbling around 
its axis of rotation.
This component of the Earth's orbit is called \emph{precession};
its period varies from 19 to~23\,Kyr.

\begin{figure}[ht]
\begin{center}
\resizebox{0.9\textwidth}{!}{\includegraphics{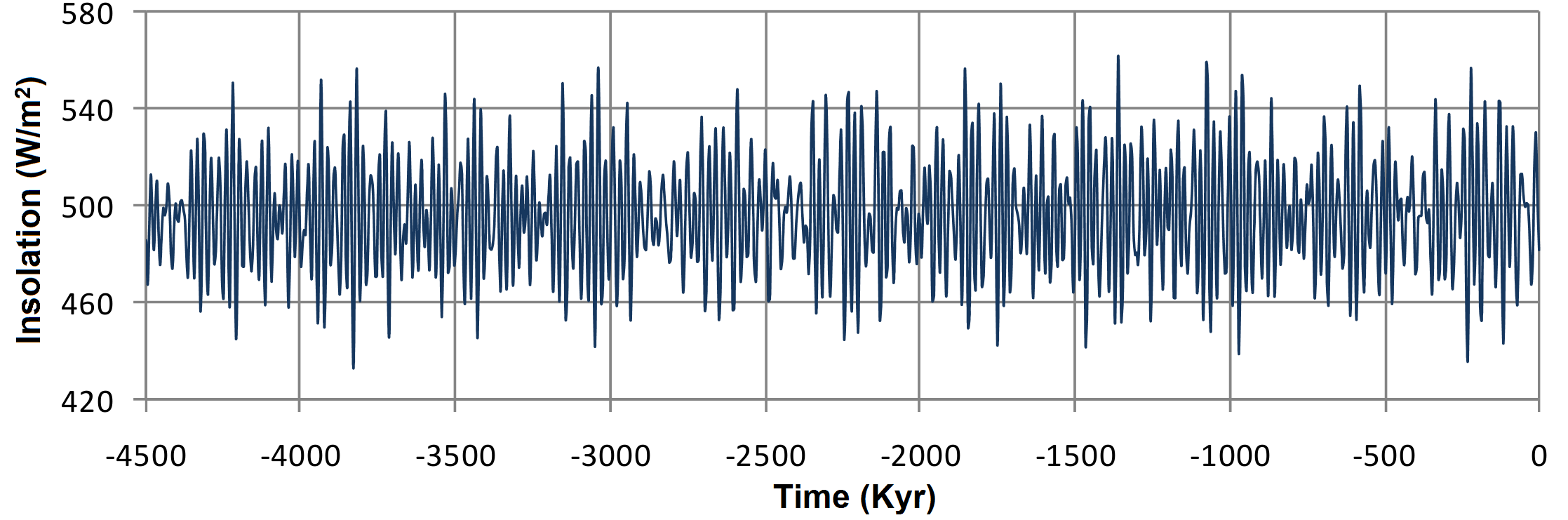}}
\caption{Time series of $Q^{65}$ during the month of July for the past 4.5 million 
years. [Data from~\cite{laskar2011la2010}]
\label{f-Q65TimeSeries}}
\end{center}
\end{figure}

Given the three orbital parameters (eccentricity,
obliquity, precession), one can compute the insolation
at any latitude and at any time of the year.
An example is given in Figure~\ref{f-Q65TimeSeries},
which shows the time series of $Q^{65}$---the average
insolation at $65^{\mathrm{o}}$\,North---during the month
of July for the past 4.5 million years;
other months show a similar behavior.
A cycle with a period of approximately 400\,Kyr is clearly visible. 
A spectral analysis reveals a dominant frequency around 21\,Kyr
coming from two clustered spikes in the power spectrum
and another, smaller frequency component at approximately~41\,Kyr.

\subsection{Atmospheric Carbon Dioxide\label{ss-CO2}}
The Pleistocene climate and, in particular, the mid-Pleistocene transition
are topics of great interest in the geosciences community.
The 41\,Kyr glacial cycles of the early Pleistocene are
commonly attributed to the 41\,Kyr cycle of Earth's 
obliquity; see, for example,~\cite{Huybers2006, Raymo2003}.
In contrast, there is less agreement on the origin of the 100\,Kyr cycles
of the late Pleistocene.

Some authors~\cite{Ghil1994, HayesImbrieShackleton1976, Imbrie1992}
attribute the 100\,Kyr cycles to the eccentricity of Earth's orbit.
However, simple energy balance considerations imply 
that variations in eccentricity are too weak to explain 
the surface temperature variations that are observed
in the paleoclimate record. 

\begin{figure}[ht]
\begin{center}
\resizebox{0.9\textwidth}{!}{\includegraphics{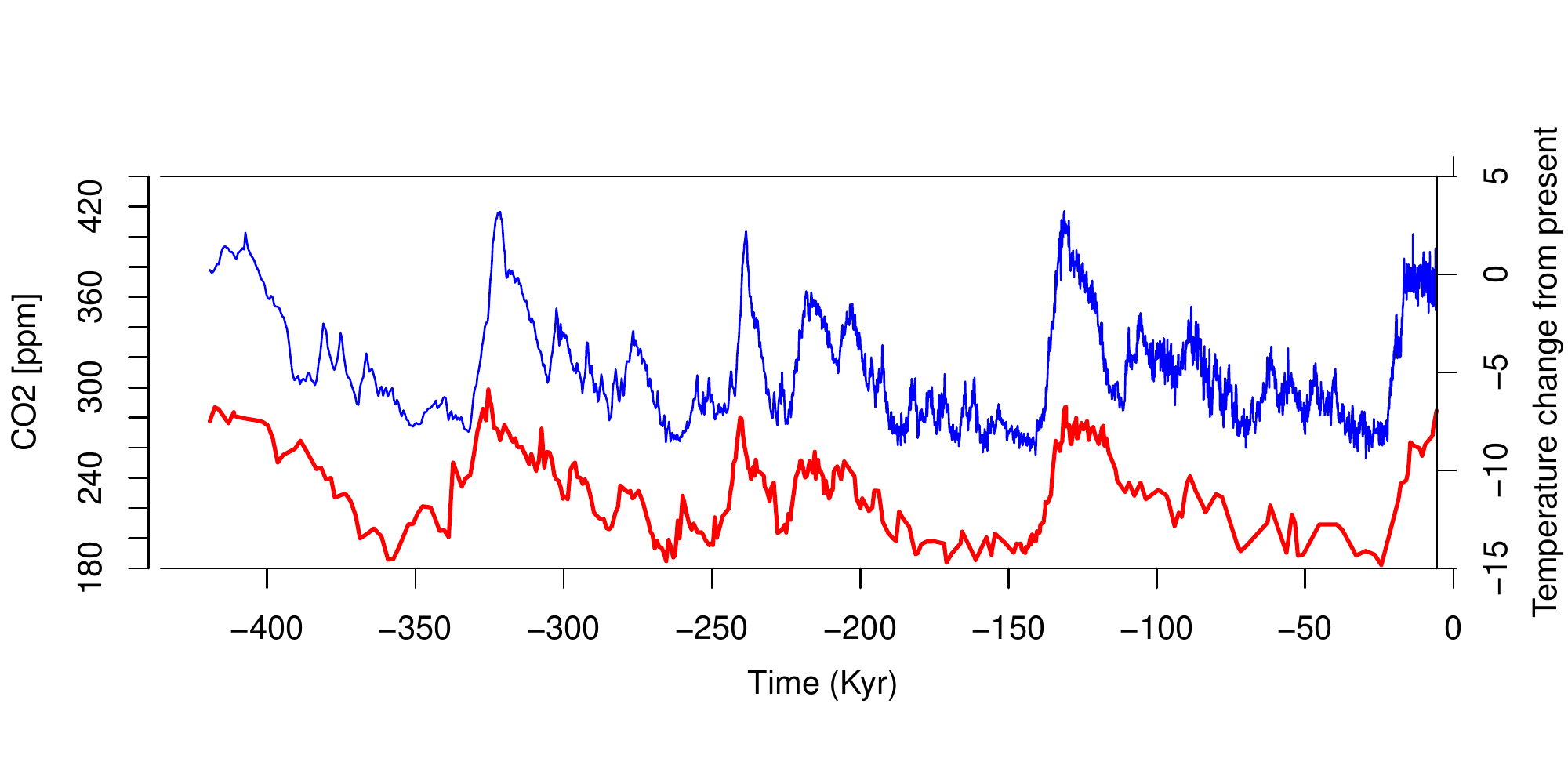}}
\caption{Correlation of global mean temperature (blue) and atmospheric \CO\ concentration (red)
for the past 420\,Kyr.
[Data from Carbon Dioxide Information Analysis Center (CDIAC) 
at Oak Ridge National Laboratory]
\label{f-CO2TemperatureCorrelation}
}
\end{center}
\end{figure}

A possible way for orbital effects to influence the Earth's surface temperature
is suggested by the greenhouse effect and the almost perfect correlation
between fluctuations in the atmospheric \CO\ concentration
and the surface temperature that is observed, for example, 
in the Vostok ice core data~\cite{Petit1999};
see Figure~\ref{f-CO2TemperatureCorrelation}.

Orbital variations have only a very weak effect
on the composition of a planet's atmosphere.
Moreover, they would affect not only \CO, 
but other atmospheric components as well,
and their effect would be negligible on the time scale
of the glacial cycles.
Hence, it is unlikely that orbital variations alone are the direct cause
of the fluctuations in atmospheric \CO\ concentrations,
so one would need some feedback mechanism that reinforces
the influence of orbital forcing on surface temperature
through changes in atmospheric \CO~\cite{rutherford2000early}. 
Changes in the carbon cycle and the climate system would then amplify each other
to produce the glacial cycles, and atmospheric \CO\ would have to play
a central role in such a feedback mechanism.

Saltzman was the first to propose a conceptual climate model
that highlighted the role of atmospheric \CO\
in the dynamics of glacial cycles~\cite{Saltzman1987}.
The model was further developed in joint work with Maasch
in a series of articles~\cite{MaaschSaltzman1990, 
SaltzmanMaasch1988, SaltzmanMaasch1991}.
In this article, we focus on the model proposed by
Maasch and Saltzman in~\cite{MaaschSaltzman1990}.

\subsection{Other Models\label{ss-OtherModels}}
There certainly is no unique way to explain the phenomenon
of glacial cycles and the Pleistocene climate in a comprehensive manner.
Other conceptual models can be found, for example, in~\cite{Ashkenazy2004, 
Ashwin2015, Clark1999, Ghil1994, Huybers2005, 
Paillard1995, Paillard1998, Paillard2001, PaillardParrenin2004,
Raymo1997, Shackleton2000, Tziperman2003}.
An interesting case was made by Huybers~\cite{Huybers2005},
who argued that the reconstruction of the temperature record
from proxy data presented in Figure~\ref{f-TempRecord5Myr}
relies on orbital assumptions and is therefore subject to bias.
Huybers developed an unbiased age model which does not rely
on orbital assumptions and showed that the late Pleistocene
glacial terminations are paced by changes in Earth's obliquity~\cite{Huybers2007}.
This theory would imply that the entire Pleistocene climate regime
can be explained by obliquity alone.
We refer the reader to~\cite{Crucifix2012} for a summary
of the current state of the art.

As a final note, we caution that any model is a mathematical construct,
and any phenomenon that results from its analysis is merely a manifestation of 
the assumptions underlying the model.
The question whether the model reflects the true cause(s) of the glacial cycles
lies outside the domain of mathematics.

\section{The Maasch and Saltzman Model \label{s-MSModel}}
\setcounter{equation}{0}
\setcounter{table}{0}
\setcounter{figure}{0}
The Pleistocene climate model proposed by
Maasch and Saltzman~\cite{MaaschSaltzman1990} 
involves five \emph{state variables}.
They are, with their associated units in square brackets:
\begin{itemize}
\item the total global ice mass, $I$ [kg],
\item the atmospheric \CO\ concentration, $\mu$ [ppm],
\item the volume of the North Atlantic Deep Water (NADW), $N$ [m${}^3$],
\item the global mean sea surface temperature (SST), $\tau$ [K], and
\item the mean volume of permanent (summer) sea ice, $\eta$ [m${}^3$].
\end{itemize}
The volume of the NADW is a measure of the strength
of the global oceanic circulation (\emph{thermohaline circulation}, THC).
We can think of $N$ as a measure of the strength of the oceanic \CO\ pump,
since the oceanic \CO\ pump is an integral part of the THC.
The other variables are self-explanatory.

The state variables vary with time, 
albeit on rather different time scales.
The total global ice mass,
atmospheric \CO\ concentration, and
NADW vary on the order of thousands of years,
while the SST and summer sea ice
vary on the order of decades or centuries.
Here the focus is on the slow time scale,
where we assume that the fast variables
equilibrate essentially instantaneously.
That is, the long-term dynamics of the climate system
are described in terms of $I$, $\mu$, and $N$
(the \emph{prognostic variables});
$\tau$ and $\eta$ are \emph{diagnostic variables},
which follow the prognostic variables in time.

\subsection{Model Formulation\label{ss-ModelFormulation}}
The climate model is formulated in terms of 
\emph{anomalies}---deviations from long-term averages, 
which are indicated with a prime.
The governing equations follow from plausible physical principles,
which are detailed in~\cite[\S 2]{SaltzmanMaasch1988}.

The global mean SST ($\tau$) and the mean volume
of permanent sea ice ($\eta$)
vary with the total global ice mass ($I$)
and the atmospheric \CO\ concentration ($\mu$)
but are independent of the NADW ($N$);
in particular, $\tau$ decreases as $I$ increases  or $\mu$ decreases,
while $\eta$ increases as $I$ increases or $\mu$ decreases.
To leading order, the dependences are linear, so
\begin{equation}
\begin{split}
\tau' &= - \alpha I' + \beta \mu' , \\
\eta' &= e_I I '- e_\mu \mu' .
\end{split}
\label{EoS}
\end{equation}
In the absence of external forces, the governing equations
for $I'$, $\mu'$, and $N'$ are
\begin{equation}
\begin{split}
\frac{dI'}{dt'} &= - s_1 \tau' - s_2 \mu' + s_3 \eta' - s_4 I'  , \\
\frac{d\mu'}{dt'} &= r_1 \tau' - r_2 \eta' - (r_3 - b_3 N') N'- (r_4 + b_4 N'^2) \mu' - r_5 I' , \\
\frac{dN'}{dt'} &= - c_0 I' - c_2 N' .
\end{split}
\label{ODE1}
\end{equation}
Time $t'$ is measured in units of one year [yr].
The coefficients in these equations are positive (or zero).
Maasch and Saltzman also included an external forcing term
related to $Q^{65}$ in the equation for $I'$, but since we are interested 
in the internal dynamics of the system, we don't include external forcing.

The physical assumptions underlying these equations are:

$\bullet$~The prognostic variables relax
to their respective long-term averages,
so their anomalies tend to zero as time increases;
in particular, $I'$ and $N'$ decay at a constant rate, 
while the decay rate of $\mu'$ increases quadratically with $N'$
~\cite[\S 2.V]{SaltzmanMaasch1988}.

$\bullet$~If the SST exceeds its mean value ($\tau' > 0$),
the total global ice mass decreases
and the atmospheric \CO\ concentration
increases (due to outgassing);
if the SST is less than its mean value ($\tau' < 0$),
the opposite happens.
The coupling is linear to leading order.

$\bullet$~Since \CO\ is a greenhouse gas,
an increase in the atmospheric \CO\ concentration
leads to a warmer climate and thus a decrease in the
total global ice mass.

$\bullet$~If the volume of permanent sea ice
exceeds its mean value ($\eta' > 0$),
the total global ice mass increases and
the atmospheric \CO\ concentration decreases;
if the volume of permanent sea ice is less
than its mean value ($\eta' < 0$),
the opposite effect happens.
The coupling is linear to leading order.

$\bullet$~A greater-than-average total global ice mass ($I' > 0$)
negatively affects both the atmospheric \CO\ concentration
and the strength of the North Atlantic overturning circulation;
a less-than-average total global ice mass ($I' < 0$)
has the opposite effect.
The coupling is linear to leading order.

$\bullet$~The atmospheric \CO\ concentration decreases
as the strength of the North Atlantic overturning circulation increases,
but the coupling weakens (strengthens) as the strength
of the NADW is above (below) average~\cite[\S 2.III\,a,b]{SaltzmanMaasch1988}.
\medskip

Upon substitution of the expressions~\eqref{EoS},
the governing equations~\eqref{ODE1} become
\begin{equation}
\begin{split}
\frac{dI'}{dt'} &= - a_0 I' - a_1 \mu' , \\
\frac{d\mu'}{dt'} &= - b_0 I' + (b_1 - b_4 N'^2) \mu' - (b_2 - b_3 N') N' , \\
\frac{dN'}{dt'} &= - c_0 I' - c_2 N' ,
\end{split}
\label{ODE2}
\end{equation}
where
\begin{equation*}
\begin{aligned}
a_0 &= s_4 - (\alpha s_1 + e_I s_3),  &
a_1 &= s_2 + \beta s_1 + e_\mu s_3,  & \\
b_0 &=  r_5 + \alpha r_1 + e_I r_2, &
b_1 &=  \beta r_1 +  e_\mu r_2 - r_4, & 
b_2 &= r_3 .
\end{aligned}
\end{equation*}
Following~\cite{MaaschSaltzman1990, SaltzmanMaasch1988},
we take $b_0 = 0$ 
and assume that the remaining coefficients
are all positive.
Note that $a_0$ and $b_1$ involve positive
as well as negative contributions, so the implicit
assumption is that the positive contributions dominate.

\subsection{Nondimensional Model\label{ss-NondimensionalModel}}
Next, we reformulate the system of equations~\eqref{ODE2}
by rescaling time, 
\begin{equation}
t = a_0 t',
\label{rescaleTime}
\end{equation}
and introduce dimensionless variables,
\begin{equation}
X = \frac{I'}{\hat{I}} , \quad Y = \frac{\mu'}{\hat{\mu}} , \quad Z = \frac{N'}{\hat{N}} ,
\label{XYZ}
\end{equation} 
where $\hat{I}$, $\hat{\mu}$, and $\hat{N}$ are reference values
of $I$, $\mu$, and $N$, respectively.
Since $a_0 \approx 1.00 \cdot 10^{-4} \, \mathrm{yr}^{-1}$,
a unit of $t$ corresponds to (approximately) 10\,Kyr.

The governing equations for $X$, $Y$, and $Z$ are
\begin{equation}
\begin{split}
\dot{X} &= - X - \hat{a}_1 Y , \\
\dot{Y} &= (\hat{b}_1 - \hat{b}_4 Z^2) Y - (\hat{b}_2 - \hat{b}_3 Z) Z , \\
\dot{Z} &= - \hat{c}_0 X - \hat{c}_2 Z ,
\end{split}
\label{MS-XYZ}
\end{equation}
where the dot $\dot{ }$ indicates differentiation with respect to~$t$.
Recall that we have set $b_0 = 0$.
The remaining coefficients are dimensionless combinations
of the physical parameters in the system of equations~\eqref{ODE2},
\[
\hat{a}_1 = \frac{a_1 \hat{\mu}}{a_0 \hat{I}} , \;
\hat{b}_1 = \frac{b_1}{a_0} , \;
\hat{b}_2 = \frac{b_2 \hat{N}}{a_0 \hat{\mu}} , \;
\hat{b}_3 = \frac{b_3 \hat{N}^2}{a_0 \hat{\mu}} , \;
\hat{b}_4 = \frac{b_4 \hat{N}^2}{a_0} , \;
\hat{c}_0 = \frac{c_0 \hat{I}}{a_0 \hat{N}} , \;
\hat{c}_2 = \frac {c_2}{a_0} .
\]
A rescaling of the variables $X$, $Y$, and $Z$,
\begin{equation}
x = \left( (\hat{c}_0 / \hat{c}_2) \surd{\hat{b}_4} \right) X , \quad
y = \left( \hat{a}_1 (\hat{c}_0 / \hat{c}_2) \surd{\hat{b}_4} \right) Y, \quad
z = \left(\surd{\hat{b}_4} \right) Z ,
\label{xyz}
\end{equation}
leads to the following dynamical system for the triple $(x, y, z)$:
\begin{equation}
\begin{split}
\dot{x} &= - x - y , \\
\dot{y} &= ry - pz + s z^2 - y z^2 , \\
\dot{z} &= - qx -q z .
\end{split}
\label{MS-xyz}
\end{equation}
The coefficients $p$, $q$, $r$, and $s$ are combinations of the physical parameters,
\begin{equation}
p =  \frac{\hat{a}_1 \hat{b}_2 \hat{c}_0}{\hat{c}_2} , \quad
q = \hat{c}_2, \quad
r = \hat{b}_1 , \quad
s = \frac{\hat{a}_1 \hat{b}_3 \hat{c}_0}{\hat{c}_2 \surd{\hat{b}_4}} .
\end{equation}
The coefficients are assumed to be positive, with $q > 1$.
The system of equations~\eqref{MS-xyz} is the model
proposed by Maasch and Saltzman in~\cite[Eqs.~(4)--(6)]{MaaschSaltzman1990}.
Note that the model considered here does not include external forcing, so it
describes the \emph{internal dynamics} of the climate system.

\subsection{Discussion\label{ss-Discussion}}
The system of equations~\eqref{MS-xyz} is
what is known as a \emph{conceptual model}.
Its derivation involves physical arguments,
but there is no guarantee that it corresponds to what
actually happened in the climate system during the Pleistocene.
Its sole purpose is to describe a possible mechanism
that explains the observed behavior of the glacial cycles.

Loosely speaking, we identify $x$, $y$, and $z$
with the anomalies of the total amount of ice,
the atmospheric \CO\ concentration,
and the volume of the NADW 
(the strength of the oceanic \CO\ pump),
respectively.
Time is normalized and expressed in units of
the characteristic time of the total global ice mass,
typically of the order of 10\,Kyr.

Because of the various transformations needed to get from
the physical system~\eqref{ODE2} to the dynamical system~\eqref{MS-xyz},
it is difficult to relate the parameters to actual physical processes.
The best we can do is look at their effect on the possible solutions.
For example, a nonzero value of the parameter $s$ renders the problem asymmetric,
so $s$ is introduced to achieve the observed asymmetry of the glacial cycles.
The coefficient $q$ is the characteristic time of NADW
(expressed in units of the characteristic time of the total
global ice mass).
The assumption $q > 1$ implies that NADW changes
on a faster time scale than the total global ice mass,
and as $q$ increases, this change occurs on an increasingly faster time scale.
If we rewrite the second equation as $\dot{y} = (r - z^2) y - p z - sz^2$,
we see that the growth rate $r$ of the atmospheric \CO\ concentration
is balanced by the anomaly of NADW.
Lastly, the coefficient $p$ expresses the sensitivity
of the atmospheric \CO\ concentration to NADW.

Conceptually, the following sequence of events hints at
the possible existence of periodic solutions:
(i)~As the amount of \CO\ in the atmosphere increases
and $y$ becomes positive, the total amount of ice decreases
and $x$ becomes negative (first equation);
(ii)~As $x$ becomes negative, the volume of NADW
increases and $z$ becomes positive (third equation);
(iii)~As $z$ becomes positive, the amount of atmospheric \CO\
decreases and $y$ becomes negative (second equation).
This is the first part of the cycle.
Once $y$ is negative, the opposite effects happen.
(iv)~The total global ice mass starts to increase again
and $x$ becomes positive (first equation);
(v)~As a result, the volume of NADW decreases
and eventually $z$ becomes negative (third equation);
(vi)~Once $z$ is negative, $y$ starts to increase again (second equation). 
This completes the full cycle and sets the stage for the next cycle.
Of course, these arguments do not guarantee the existence
of a periodic cycle and do not say anything about its period.
The particulars will depend critically on the parameter values.

\subsection{Computational Results\label{ss-ComputationalResults}}
\begin{figure}[ht]
\begin{center}
\resizebox{0.8\textwidth}{!}{\includegraphics{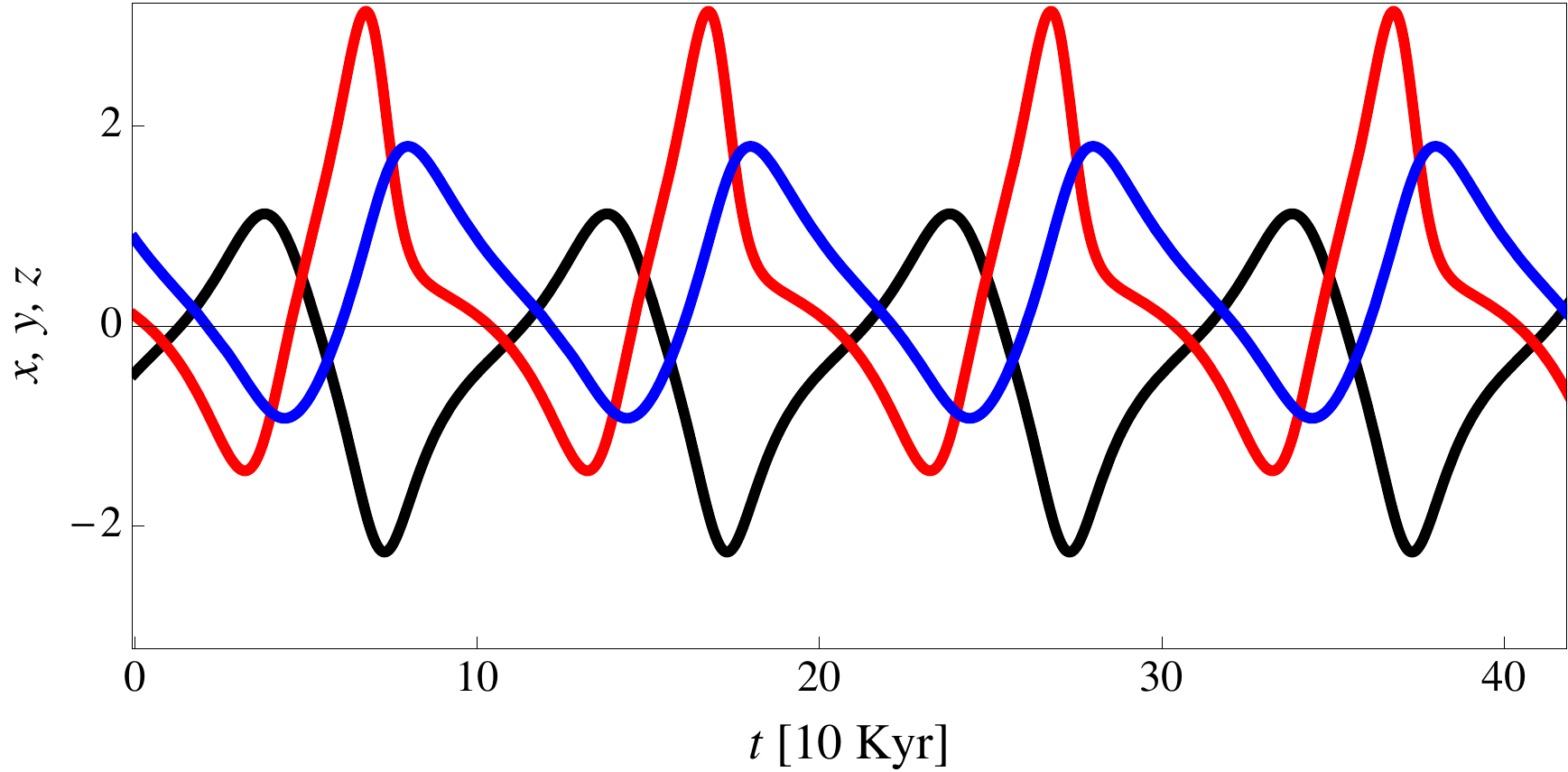}}
\caption{Limit cycle of~\eqref{MS-xyz}
at $p =1.0$, $q = 1.2$, $r = 0.8$, $s = 0.8$.
\label{f-figNoOrbital}}
\end{center}
\end{figure}

Maasch and Saltzman found computationally
that the system~\eqref{MS-xyz} generates a limit cycle
with a 100\,Kyr period at the parameter values
$p =1.0$, $q = 1.2$, $r = 0.8$, and $s = 0.8$.
The limit cycle is shown in Figure~\ref{f-figNoOrbital}.
The three curves represent the total ice mass (black), 
the atmospheric \CO\ concentration (red),
and the volume of NADW (blue)
in arbitrary units.
Each cycle is clearly asymmetric: 
a rapid deglaciation is followed by a slow glaciation.
Also, the three variables are properly correlated:
as the concentration of atmospheric \CO\  (a greenhouse gas)
increases, the climate gets warmer and the total ice mass decreases;
as the volume of NADW increases, the strength of the
North Atlantic overturning circulation increases, 
more atmospheric \CO\ is absorbed by the ocean
and, consequently, the atmospheric \CO\ concentration decreases. 

\begin{figure}[ht]
\begin{center}
\resizebox{0.7\textwidth}{!}{\includegraphics{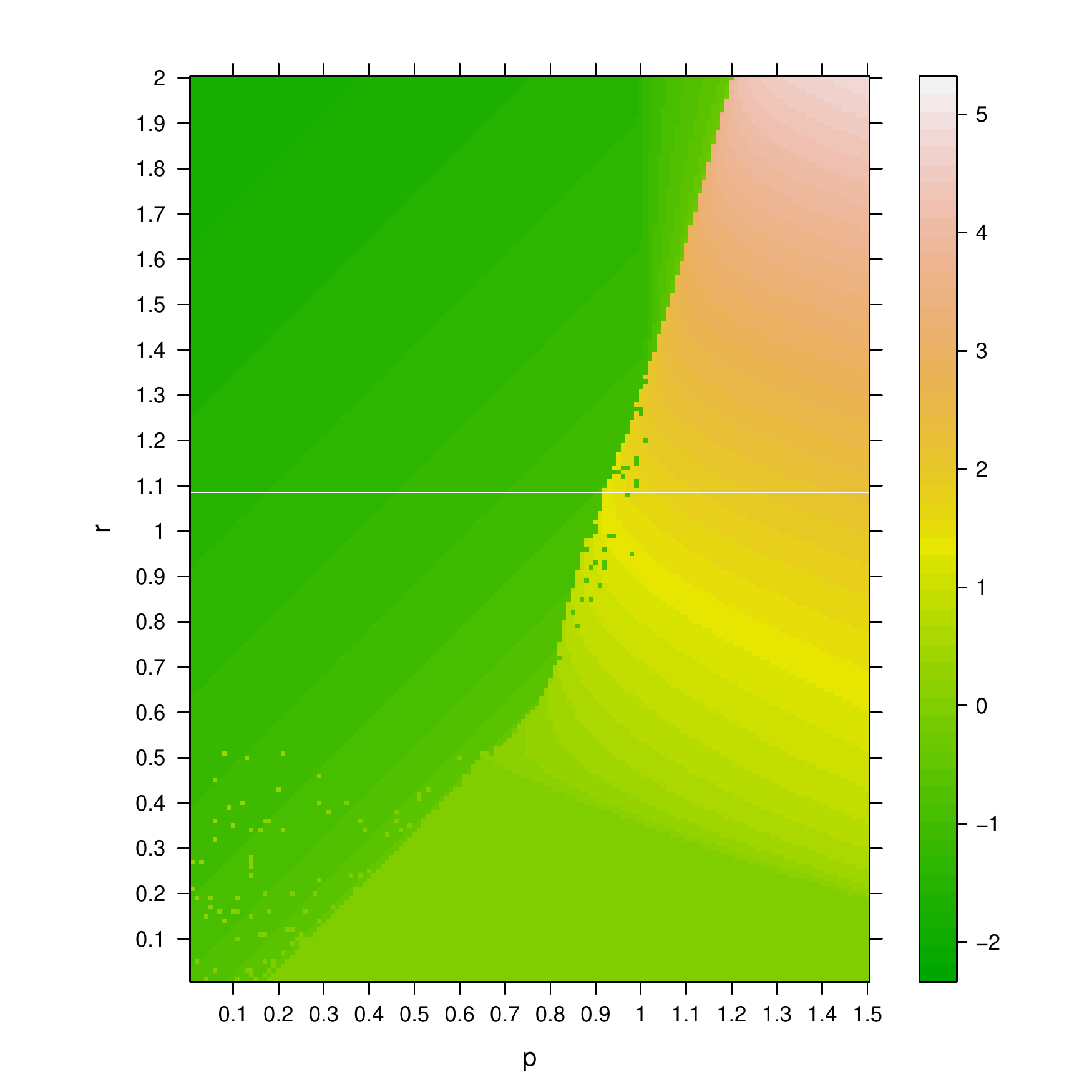}}
\caption{Color map of $\overline{x} (p,r)$ for the system~\eqref{MS-xyz} at $q = 1.2$, $s = 0.8$,
indicating convergence to an equilibrium state or a limit cycle.
\label{f-Intro-ExploreLimitCycles}}
\end{center}
\end{figure}

More detailed numerical calculations show that the Maasch-Saltzman model
possesses limit cycles in large portions of parameter space.
We integrated the system~\eqref{MS-xyz} forward in time
for a range of values of $(p, r)$, keeping $q$ and $s$ fixed
at the values $q = 1.2$ and $s = 0.8$, starting from 
a randomly chosen initial point for each $(p, r)$,
until there was a clear indication of either a limit cycle or a limit point.
We then determined the quantity $\overline{x}$ for each pair $(p, r)$,
\be
\overline{x} = \lim \sup_{t \to \infty} x(t) .
\label{xbar}
\ee
Figure~\ref{f-Intro-ExploreLimitCycles} shows the function
$(p, r) \mapsto \overline{x} (p, r)$ as a color map.
A limiting value~0 (light green) indicates convergence to the trivial state, 
a nonzero negative value (dark green) convergence to
a nontrivial equilibrium state, and a nonzero positive value
(orange or pink) convergence to a limit cycle with a finite amplitude.
Limit cycles were observed in the entire orange-colored region of the $(p, r)$ plane.
These findings, as well as other results reported by
Maasch and Saltzman in~\cite{MaaschSaltzman1990},
especially when the effects of orbital forcing are included,
suggest that the conceptual model~\eqref{MS-xyz} may indeed provide
an explanation for the Pleistocene climate record.

\section{Simplifying the Maasch--Saltzman Model\label{s-SimplifyMS}}
The system~(\ref{MS-xyz}) has four positive parameters,
$p$, $q$, $r$, and $s$, where $q > 1$. 
As noted in Section~\ref{ss-Discussion},
the parameter~$s$ introduces asymmetry into the model.
If $s=0$, the equations are invariant under reflection:
if $(x, y, z)$ is a solution, then so is $(-x, -y, -z)$.
Note furthermore that, as $q \to \infty$, the differential equation for $z$
reduces, at least formally, to the identity $z = - x$,
so the system~(\ref{MS-xyz}) becomes two-dimensional.
These observations suggest that it may be helpful to analyze
the dynamics of the Maasch--Saltzman model~(\ref{MS-xyz})
in stages, where we first focus on the special case
$q = \infty$ and $s = 0$, and then consider the effects of
finite values of $q$ and positive values of $s$.

If we set $q = \infty$ and $s = 0$, the system~(\ref{MS-xyz})
reduces formally to a two-dimensional system with $\Z_2$ symmetry,
\begin{equation}
\begin{split}
\dot{x} &= - x - y , \\
\dot{y} &= ry + px - x^2 y .
\end{split}
\label{2Dsym-xy}
\end{equation}
The state variables are $x$ and $y$, the state space is $\R^2$,
$p$ and $r$ are parameters, and the parameter space is $\R_+^2$.
Its dynamics can be analyzed rigorously and completely.

\subsection{Equilibrium States and Their Stability\label{ss-2Dsym-EquilStab}}
\begin{figure}[ht]
\begin{center}
\begin{minipage}{0.7\textwidth}
\resizebox{\textwidth}{!}{\includegraphics{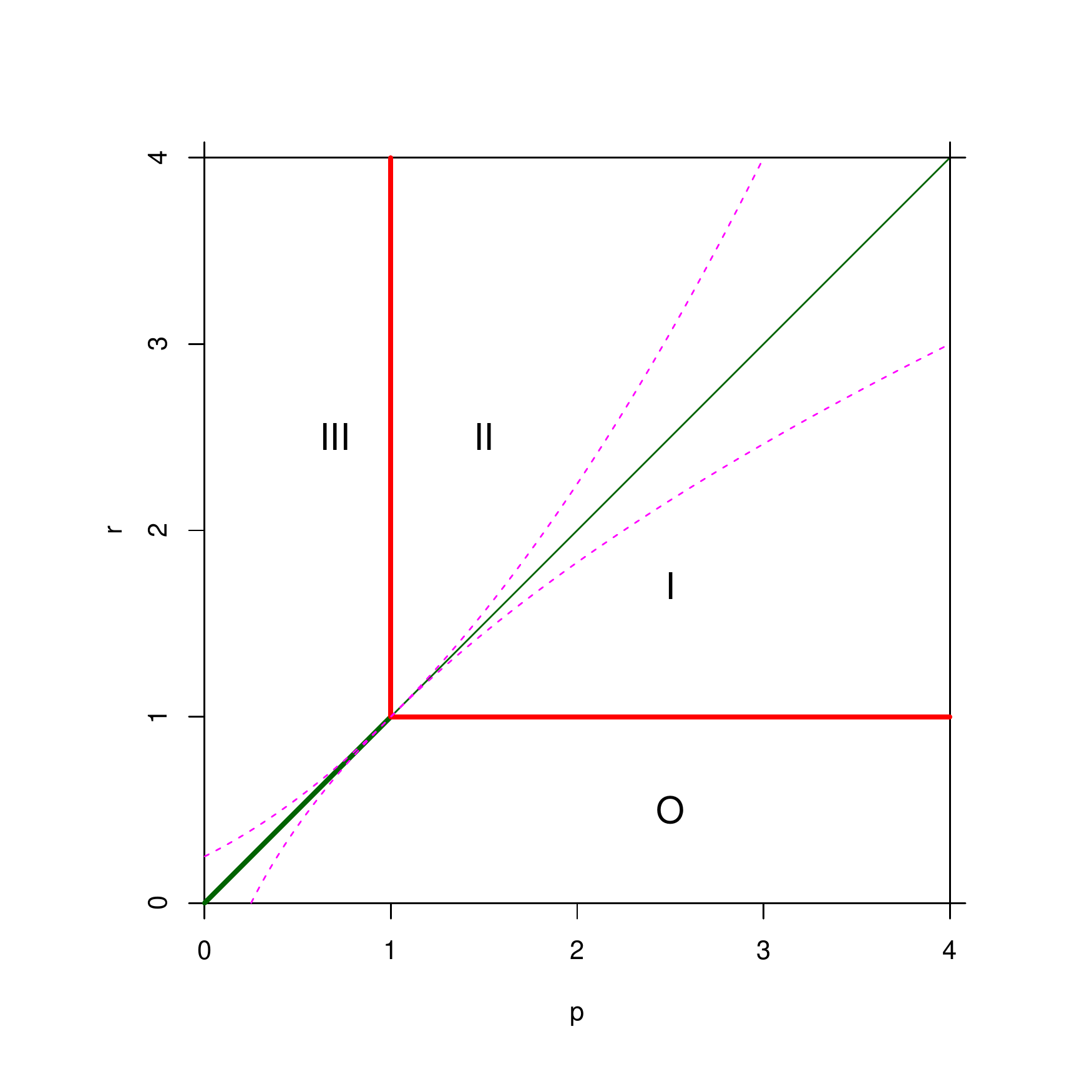}}
\end{minipage}
\begin{minipage}{0.25\textwidth}
\resizebox{0.9\textwidth}{!}{\includegraphics{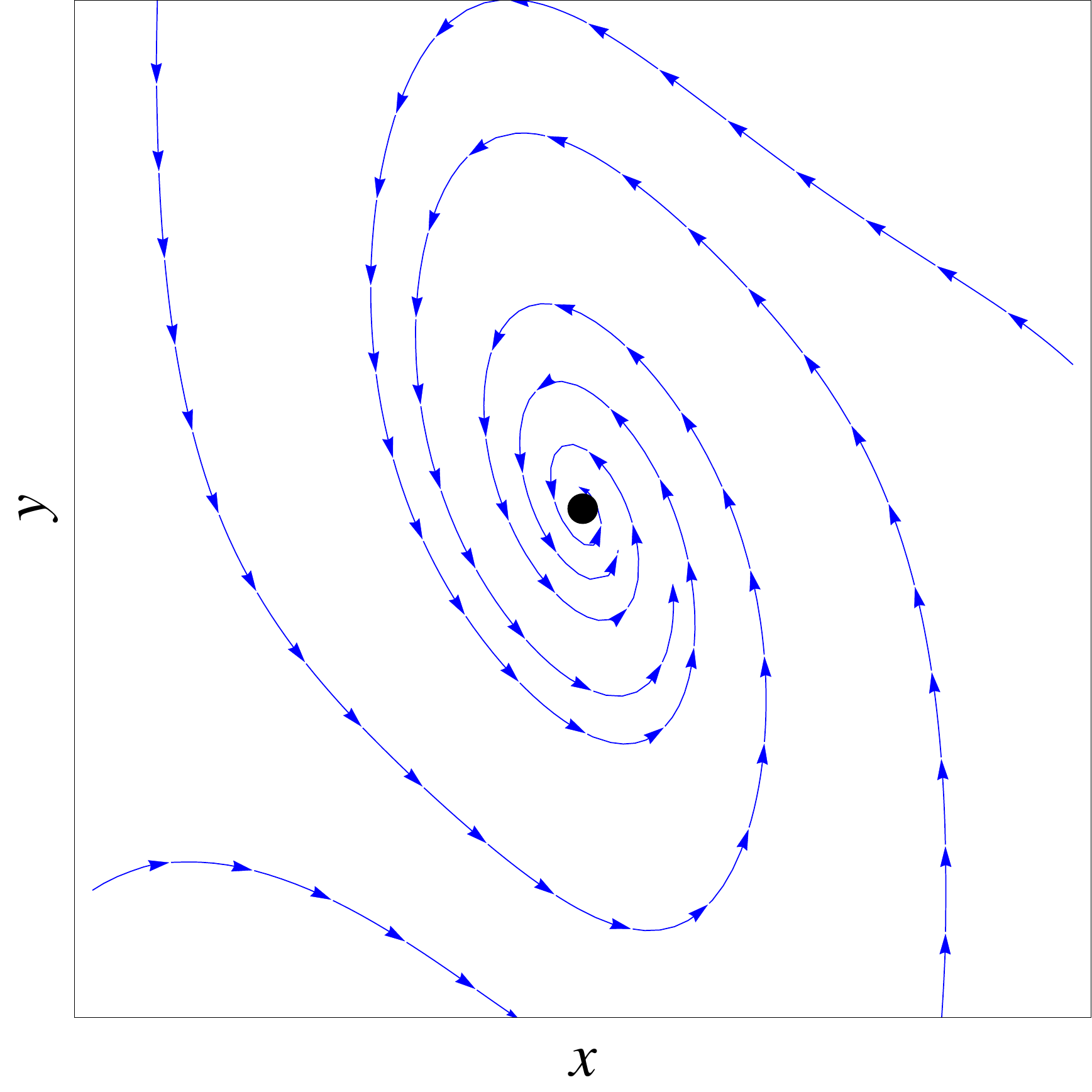}} \\
\resizebox{0.9\textwidth}{!}{\includegraphics{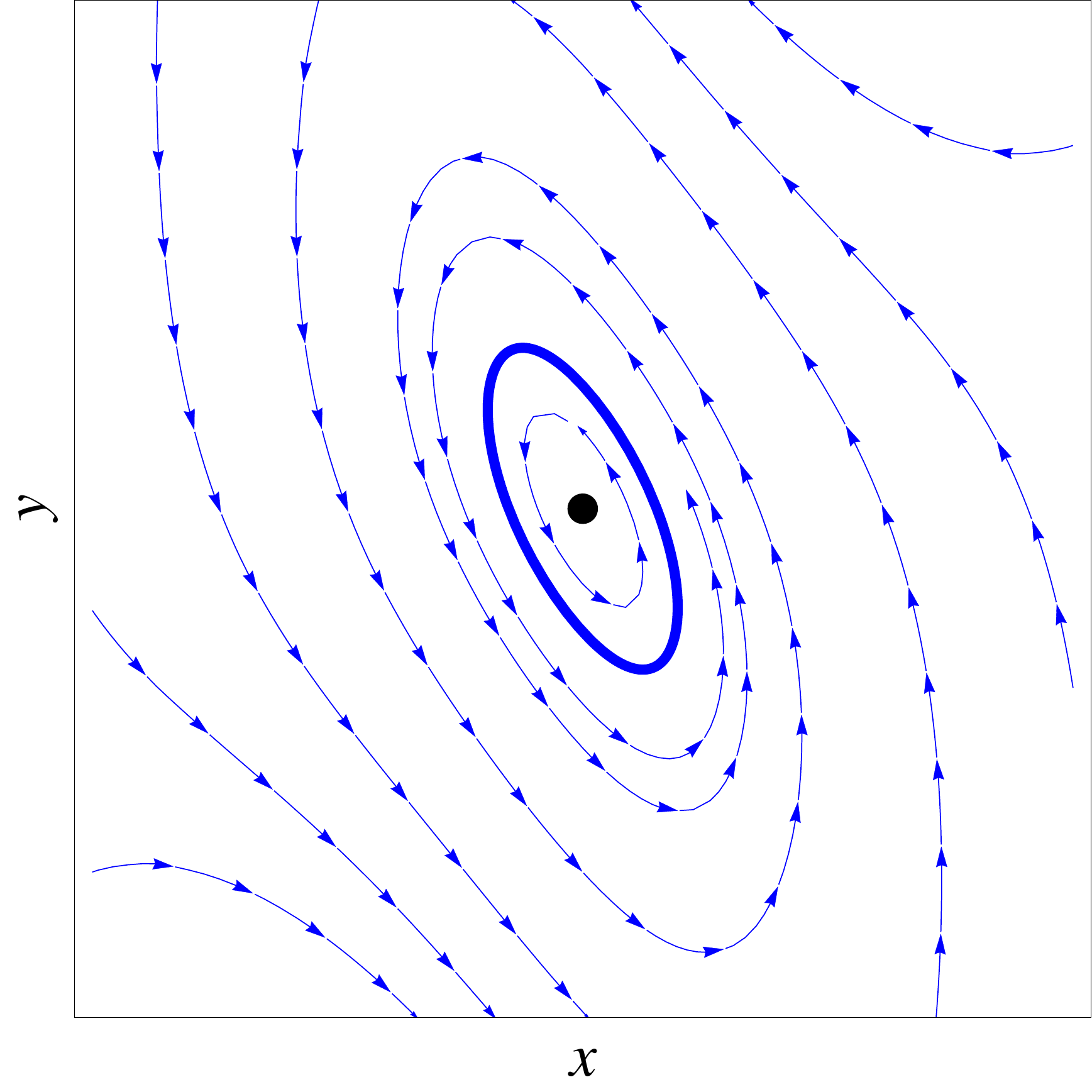}} \\
\resizebox{0.9\textwidth}{!}{\includegraphics{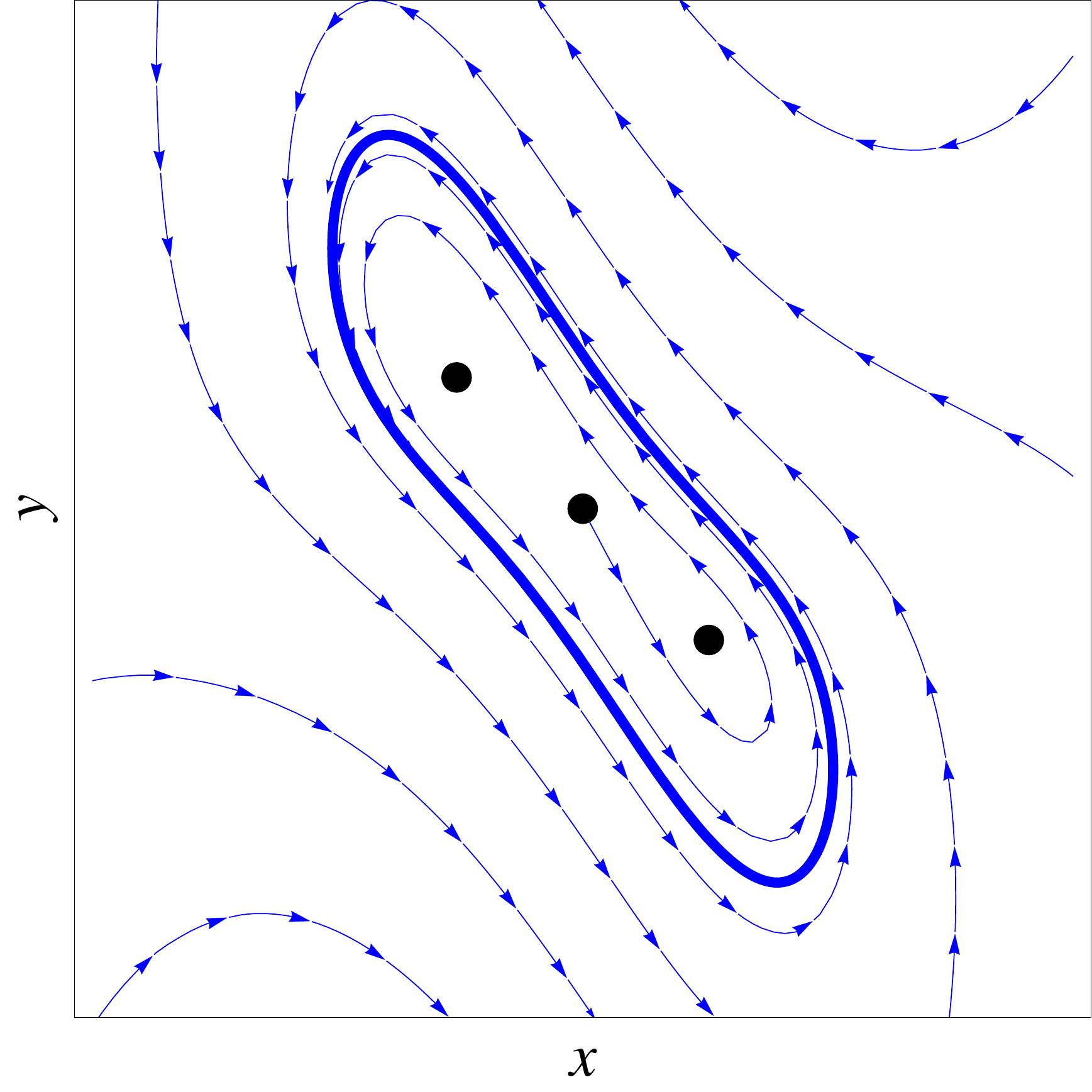}}
\end{minipage}
\end{center}
\caption{
(Left) Stability regions of the equilibrium states $P_0$, $P_1$, and $P_2$
of~(\ref{2Dsym-xy}).
(Right) Representative trajectories in region~O (top), I (middle), and II (bottom).
\label{f-2Dsym-stability}}
\end{figure}

The origin $P_0 = (0,0)$ is an equilibrium state of the system~(\ref{2Dsym-xy})
for all values of $p$ and $r$.
If $r > p$, there are two additional equilibrium states,
$P_1 =  (x_1^*, -x_1^*)$ with $x_1^* = \sqrt{r-p}$,
and  $P_2 =  (x_2^*, -x_2^*)$ with $x_2^* = - \sqrt{r-p}$.
A (local) linear stability analysis shows that
$P_0$ is stable if $0 < r < \min (p, 1)$, unstable otherwise;
$P_1$ and $P_2$ are stable if $0 < p < \min(r, 1)$, unstable otherwise.
Thus, the parameter space is partitioned into four regions,
\begin{equation}
\begin{split}
\mathrm{O} &= \{ (p, r) \in \R_+^2 : r < p \mbox{ for } p < 1, \, r < 1 \mbox{ for } p > 1 \} , \\
\mathrm{I} &= \{ (p, r) \in \R_+^2 : 1 < r < p \mbox{ for }p > 1 \} , \\
\mathrm{II}& = \{ (p, r) \in \R_+^2 : 1 < p < r \mbox{ for } r > 1 \} , \\
\mathrm{III} &= \{ (p, r) \in \R_+^2 : p < r \mbox{ for } r < 1, \, p < 1 \mbox{ for } r > 1 \}  .
\end{split}
\label{2Dsym-O-III}
\end{equation}
The regions are shown in Figure~\ref{f-2Dsym-stability},
together with representative trajectories in regions~O, I, and~II.
The diagonal $r = p$ is the locus of pitchfork bifurcations,
where $P_1$ and $P_2$ are created as stable nodes
as $(p, r)$ crosses the diagonal from region~O into region~III
or as unstable nodes as $(p, r)$ crosses the diagonal
from region~I into region~II.
The parabolic curves
$C_1 = \{ p = \tfrac14 (r+1)^2 \}$ and 
$C_2 = \{ r = \tfrac14 (p+1)^2 \}$
(dashed curves shown in purple), 
which are tangent to the diagonal $r = p$ at the point~$(1,1)$,
mark the boundaries between spirals and nodes.

\subsection{Hopf Bifurcation\label{ss-2DsymHopf}}
The system~(\ref{2Dsym-xy}) is equivalent with the Li\'{e}nard equation
\begin{equation}
\ddot{x} + g (x) \dot{x} + f (x) = 0 ,
\label{2Dsym-Lienard}
\end{equation}
where $f$ and $g$ are polynomial functions,
\begin{equation}
f (x) = x^3 - (r-p) x , \quad
g (x) = x^2 - (r-1) .
\label{2Dsym-fg}
\end{equation}
A Hopf bifurcation occurs when $f (x^*) = 0$, $f' (x^*) > 0$,
and $x^*$ is a simple zero of~$g$.
The natural frequency of oscillations is $\omega^* = \sqrt{f'(x^*)}$,
and the first Lyapunov coefficient at $x^*$ is
\begin{equation}
\ell^* = - \frac{\omega^*}{8} \left. \frac{d}{dx} \frac{g' (x)}{f' (x)} \right|_{x = x^*} ,
\label{2Dsym-Hopf-ell*}
\end{equation}
see~\cite[\S 3.5]{Kuznetsov2013}.
If $\ell^* < 0$, the Hopf bifurcation is supercritical;
if $\ell^* > 0$, it is subcritical.
The sign of $\ell^*$ is the same as that of 
$f'' (x^*) g' (x^*) - f' (x^*) g'' (x^*)$,
which in the case of the polynomial functions
$f$ and $g$ given in~(\ref{2Dsym-fg})
is the same as that of $3 (x^*)^2 + r - p$.

The red line segments in Figure~\ref{f-2Dsym-stability}
are loci of Hopf bifurcations.
On the horizontal segment at $r=1$, which is associated with $P_0$,
we have $x^* = 0$.
The sign of $\ell^*$ is the same as that of $1 - p$,
which is negative, so the Hopf bifurcation is supercritical.
The natural frequency is $\omega^* = \sqrt{p-1}$.
On the vertical segment at $p=1$, which is associated with $P_1$ and $P_2$,
we have $(x^*)^2 = r-p$.
The sign of $\ell^*$ is the same as that of $4(r - 1)$,
which is positive, so the Hopf bifurcation is subcritical.
The natural frequency is $\omega^* = \sqrt{2(r-1)}$.

\subsection{Organizing Center\label{ss-2DsymBT}}
The point $(p, r) = (1, 1)$ plays a pivotal role in understanding
the complete dynamics of the system~(\ref{2Dsym-xy}).
To see why, rotate the coordinate system
by the transformation $(x, -(x+y)) \mapsto (x,y)$.
In the new coordinates, the system~(\ref{2Dsym-xy}) is
\begin{equation}
\begin{split}
\dot{x} &= y , \\
\dot{y} &= (r - p) x + (r - 1) y - x^2 y - x^3 ,
\end{split}
\label{2Dsym-xy'}
\end{equation}
and the equilibrium points are
$P_0 = (0, 0)$, $P_1 = (x_1^*, 0)$, and $P_2 = (x_2^*, 0)$.

Let $P= (x^*, 0)$ be any of the equilibrium points,
with $x^* = 0$, $x_1^*$, or $x_2^*$.
The Jacobian of the vector field at $P$ is
\begin{equation}
\begin{pmatrix} 0 & 1 \\ r - p - 3 (x^*)^2  & r - 1 - (x^*)^2 \end{pmatrix} .
\end{equation} 
The matrix has a double-zero eigenvalue at the point $(p, r) = (1, 1)$,
so the system~(\ref{2Dsym-xy'}) undergoes a Bogdanov--Takens (BT) bifurcation.
Holmes and Rand~\cite{Holmes1980} refer to such a point as an \emph{organizing center}.
Specifically, given the $\Z_2$ symmetry of the system~(\ref{2Dsym-xy'}),
the point $(p, r) = (1, 1)$ is a BT point with $\Z_2$ symmetry;
examples of such points are discussed in~\cite[Ch.~4.2]{Carr1981}
and~\cite[\S~8.4]{kuznetsov2005practical}.

The system~(\ref{2Dsym-xy'}) can be analyzed in the neighborhood
of the organizing center by the unfolding procedure outlined in 
the original papers by Bogdanov~\cite{Bogdanov1975}
and Takens~\cite[pp.~23--30]{Takens1974} (reprinted in~\cite[Chapter~1]{Broer2001})
and described in the textbooks of
Guckenheimer and Holmes~\cite[\S 7.3]{Guckenheimer1983}
and Kuznetsov~\cite[\S 8.4]{Kuznetsov2013}.
Here, we summarize the results;
the details are given in Section~\ref{ss-BTUnfolding} below.
\medskip

\begin{figure}[ht]
\begin{minipage}{0.7\textwidth}
\begin{center}
\resizebox{\textwidth}{!}{\includegraphics{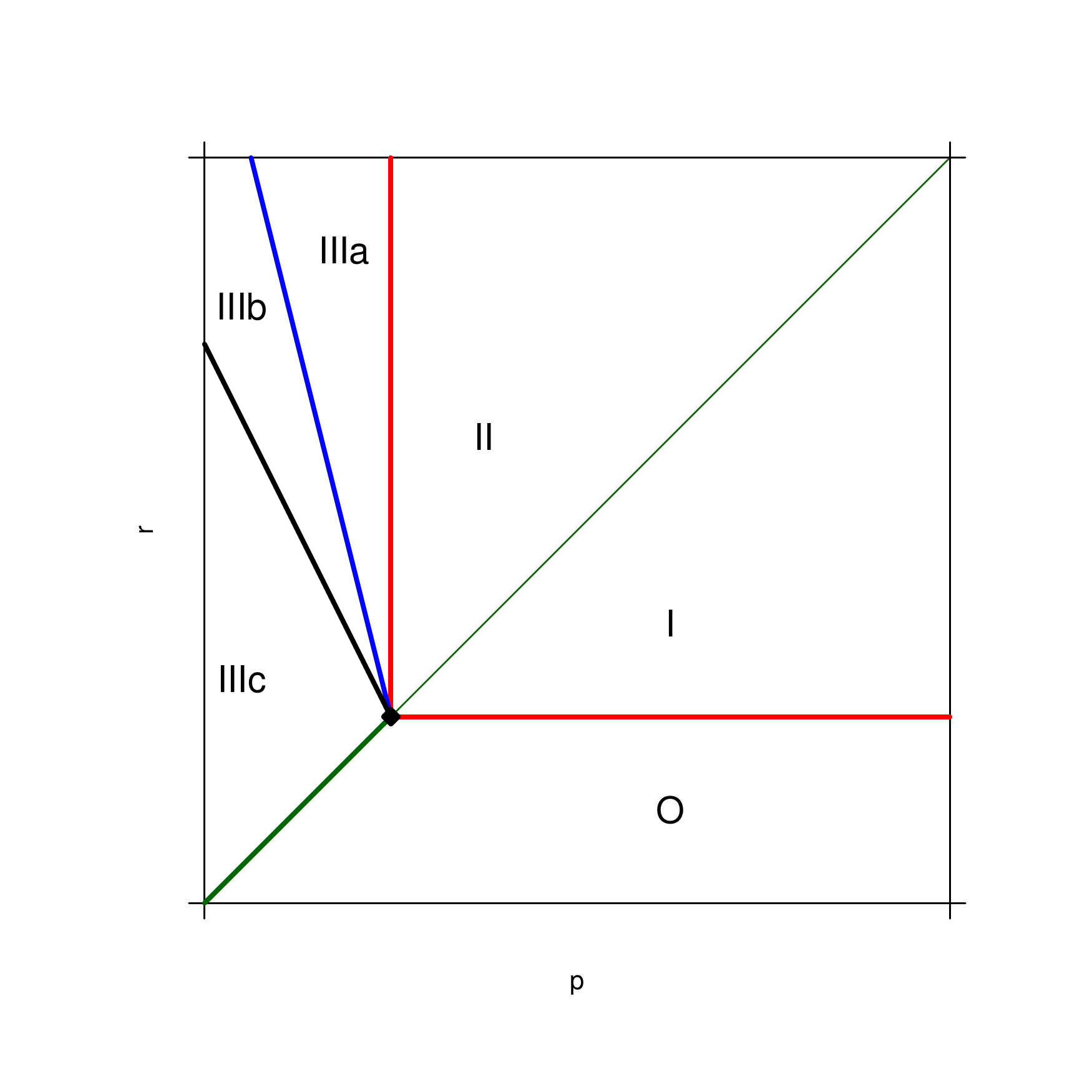}}
\end{center}
\end{minipage}
\begin{minipage}{0.25\textwidth}
\begin{center}
\resizebox{0.9\textwidth}{!}{\includegraphics{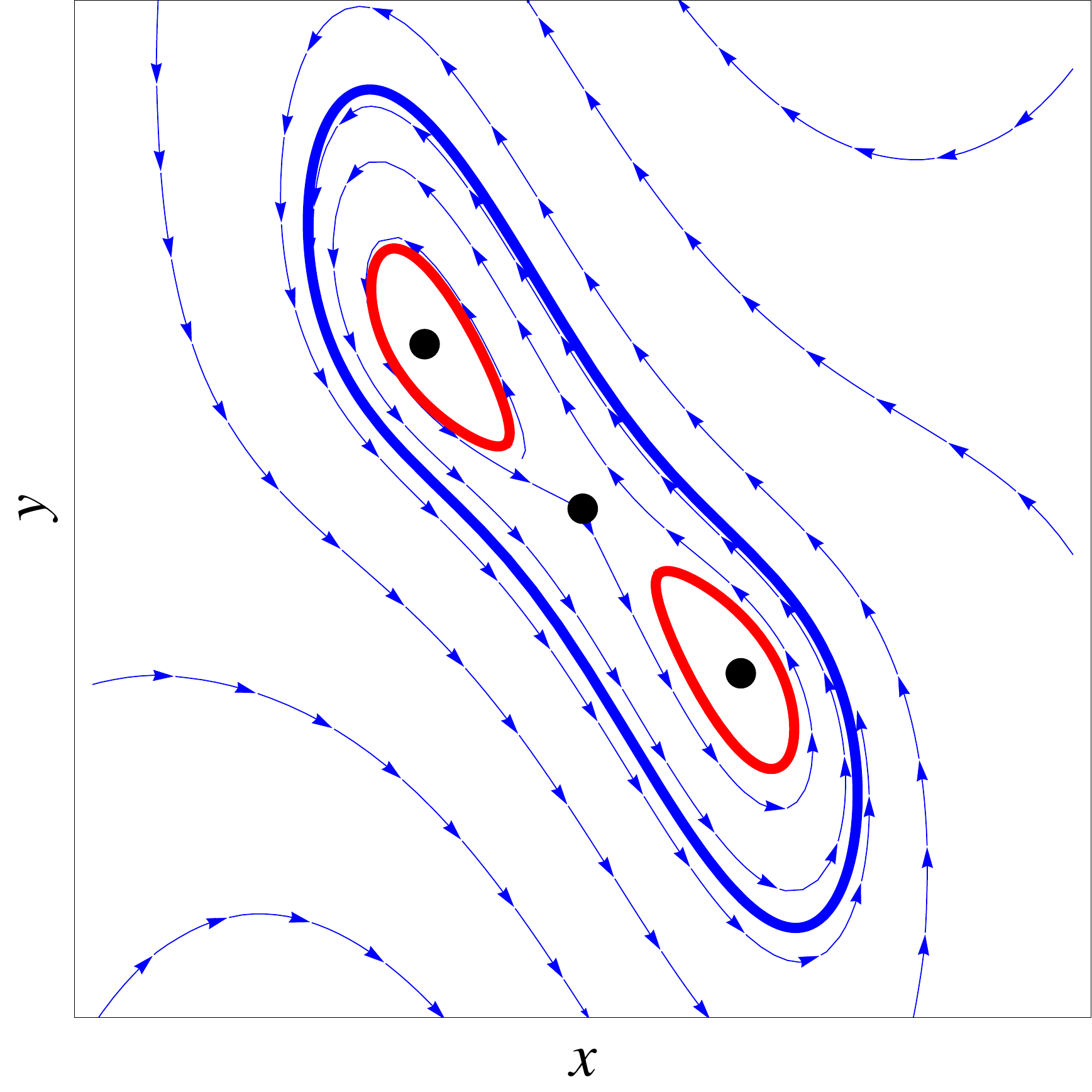}} \\
\resizebox{0.9\textwidth}{!}{\includegraphics{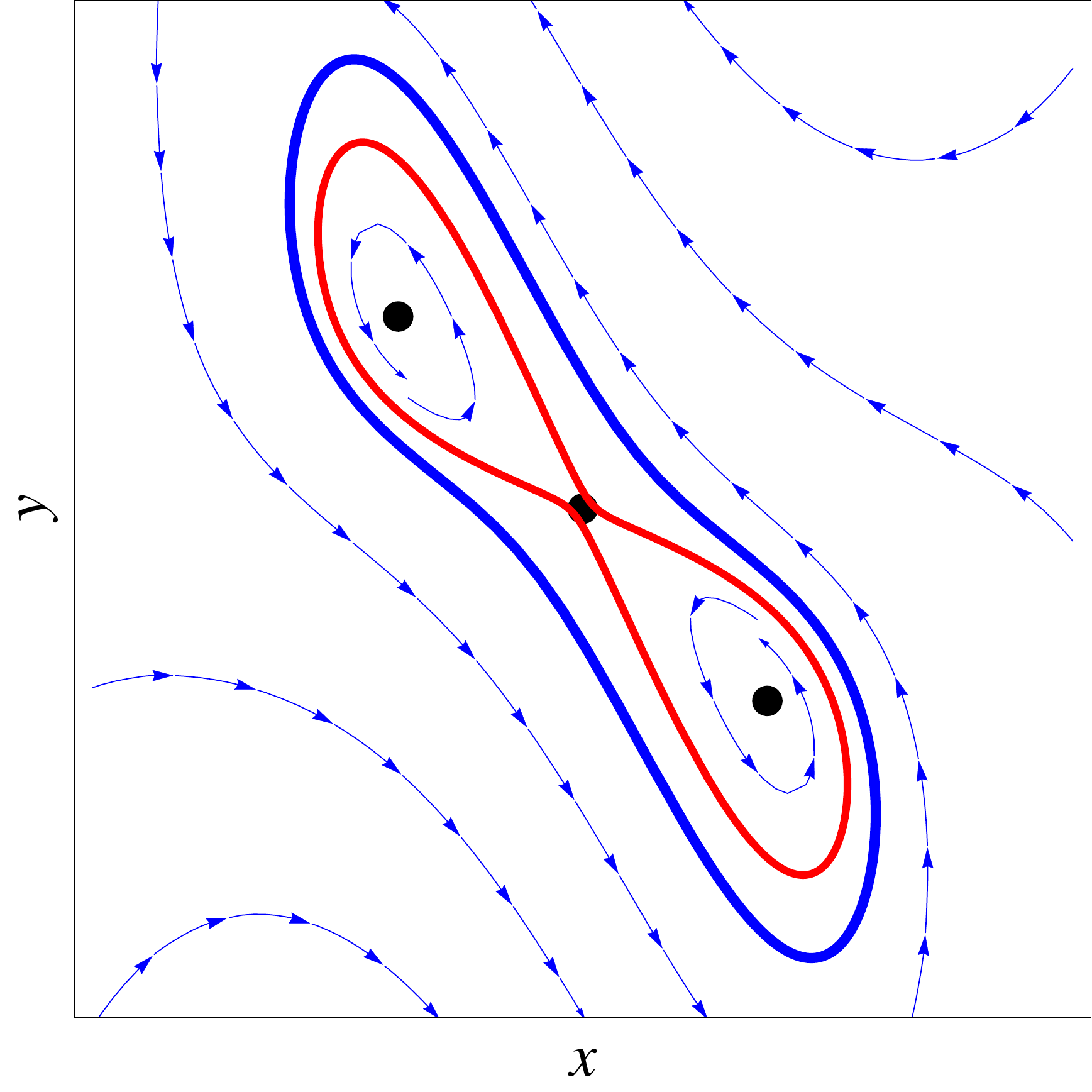}} \\
\resizebox{0.9\textwidth}{!}{\includegraphics{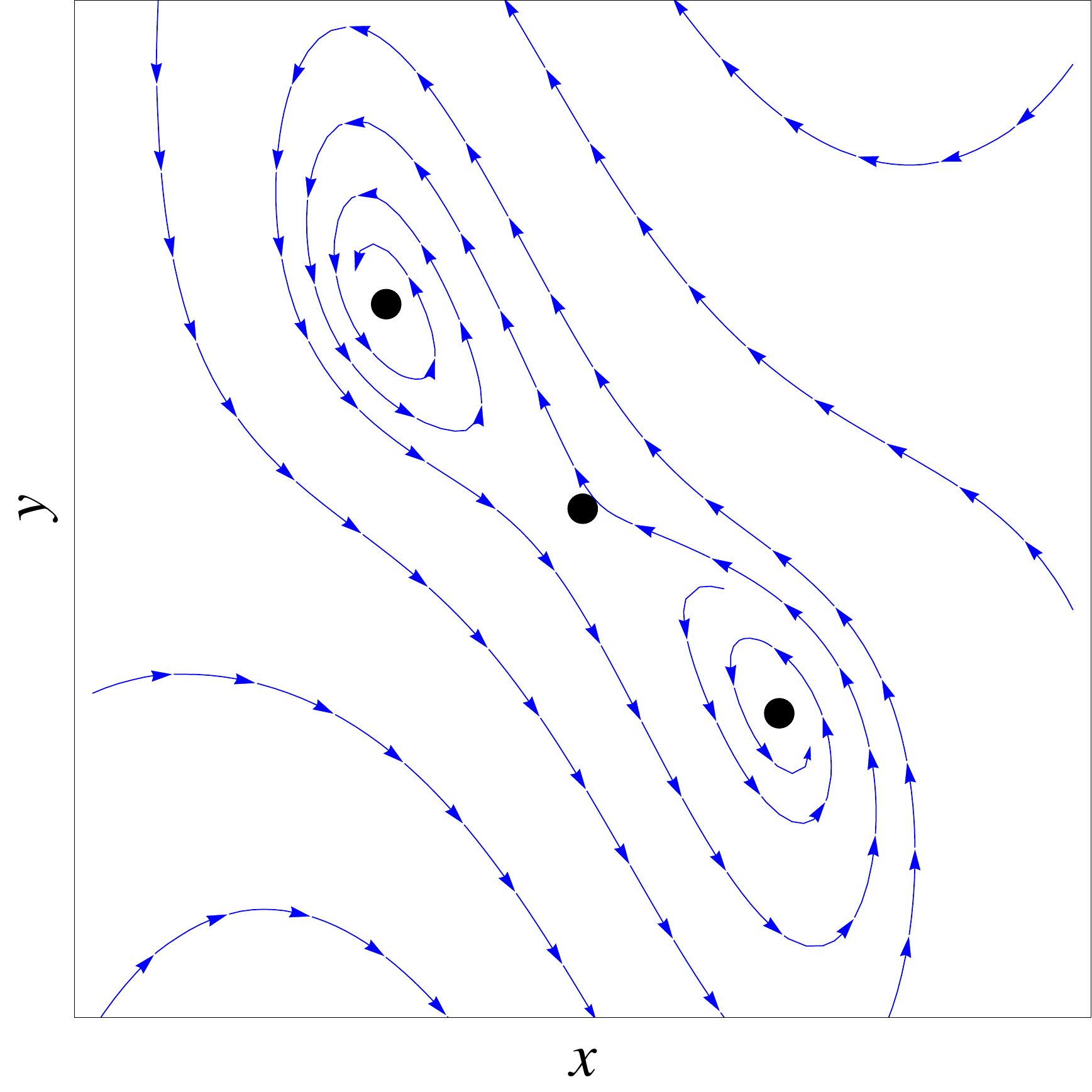}}
\end{center}
\end{minipage}
\caption{
(Left) A sketch of the bifurcation curves of~(\ref{2Dsym-xy}) near the organizing center
(the actual curves are shown in Figure~\ref{f-2Dsym-HeatPlot}).
(Right) Trajectories in region~IIIa (top), IIIb (middle), and IIIc (bottom).
\label{f-2Dsym-stability2}}
\end{figure}

Near the point $(p,r) = (1,1)$, 
region~III of Figure~\ref{f-2Dsym-stability} decomposes
into three subregions; see Figure~\ref{f-2Dsym-stability2}.
In region~IIIa, there is one stable limit cycle, 
with a pair of unstable limit cycles in its interior,
one around each of the equilibrium states
$P_1$ and~$P_2$.
As $(p, r)$ transits from region~IIIa into region~IIIb,
the two unstable periodic solutions merge to become
a pair of unstable homoclinic orbits to the saddle~$P_0$.
This homoclinic bifurcation curve (shown in blue)
is tangent to the line $r - 1 = - 4 (p - 1)$
at the organizing center $(1,1)$.
In region~IIIb, there is one stable limit cycle 
with an unstable limit cycle in its interior.
As $(p, r)$ transits from region~IIIb into region~IIIc,
there is a curve of saddle-node bifurcations of limit cycles,
along which the stable and unstable limit cycles disappear.
This curve (shown in black) is tangent to the line $r - 1 \approx - 3.03 (p - 1)$
at the organizing center $(1,1)$.
In region~IIIc, only the three equilibrium states remain,
$P_0$ as an unstable saddle, $P_1$ and~$P_2$ 
as stable spirals or nodes.

\subsection{Computational Results\label{ss-2DsymComputational}}
\begin{figure}[ht]
\begin{center}
\resizebox{0.7\textwidth}{!}{\includegraphics{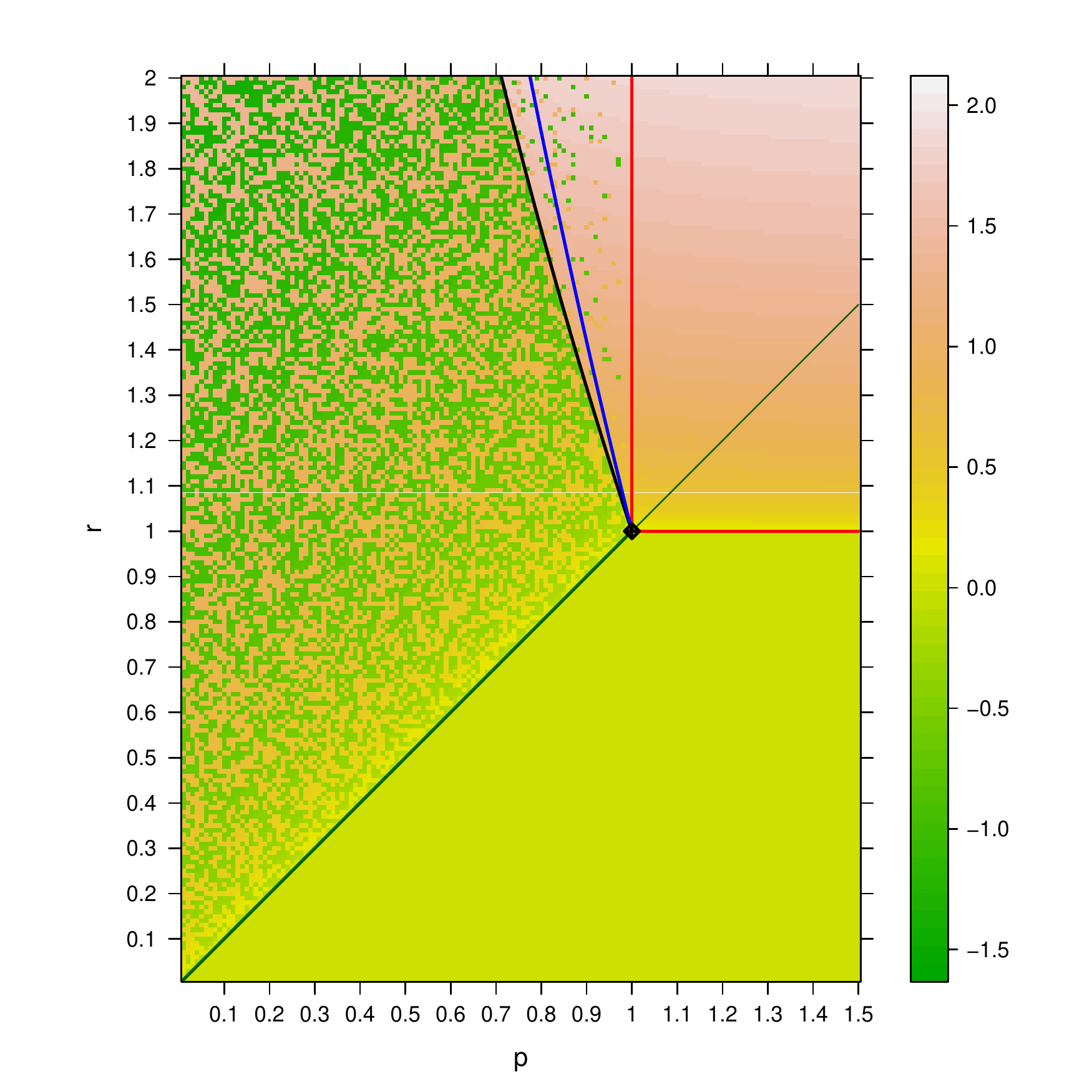}}
\caption{
Color map of $\overline{x} (p,r)$ and bifurcation curves for the system~(\ref{2Dsym-xy}).
\label{f-2Dsym-HeatPlot}}
\end{center}
\end{figure}

To complement the analysis, we performed an integration of
the system~(\ref{2Dsym-xy}) forward in time for a range of values of $(p, r)$,
following the procedure described in Section~\ref{ss-ComputationalResults}
for Figure~\ref{f-Intro-ExploreLimitCycles},
to determine the quantity~$\overline{x} = \lim \sup_{t \to \infty} x (t)$,
as in~(\ref{xbar}).
Figure~\ref{f-2Dsym-HeatPlot}  shows the function
$(p, r) \mapsto \overline{x} (p, r)$ as a color map,
together with the bifurcation curves obtained with AUTO.
A limiting value~0 (light green) indicates convergence to the trivial state, 
a nonzero negative value (dark green) convergence to
a nontrivial equilibrium state, and a nonzero positive value (orange or pink)
convergence to a limit cycle with a finite amplitude.

The stability region~O of the trivial state $P_0$ is clearly visible,
as is its Hopf bifurcation curve.
As one crosses the Hopf bifurcation curve in the direction of increasing~$r$,
the color changes slowly toward increasing values of $\overline{x}$,
indicating a supercritical Hopf bifurcation and periodic orbits
with amplitudes $\mathcal{O} (\sqrt{r-1})$.

Throughout regions~I and~II, the color map changes from orange to pink
as $r$ increases, indicating that the solutions all approach the stable limit
cycle that surrounds $P_0$; cf.~Figure~\ref{f-2Dsym-stability}.

In region~IIIa, there is a similar shift to pink as $r$ increases. 
One also sees some green and orange patches in region~IIIa, 
indicating that some of the randomly chosen initial conditions lie
in the basins of attraction of the stable equilibrium states $P_1$ and $P_2$. 
Next, in region~IIIb, the color map has largely the same characteristics
as in region~IIIa, corresponding to the fact that solutions
with initial conditions that lie inside the large unstable limit cycle
approach one of the stable equilibria (green or orange),
and those with initial conditions outside the unstable limit cycle
approach the large stable limit cycle (pink). 
Finally, in region~IIIc, the color map consists entirely of
green and orange, indicating that all of the solutions are
attracted either to $P_1$ or to $P_2$, as expected
since there are no stable limit cycles in region~IIIc.

\subsection{Bogdanov--Takens Unfolding\label{ss-BTUnfolding}}
In Section~\ref{ss-2DsymBT}, we presented the results
of a bifurcation analysis of the system~(\ref{2Dsym-xy'})
in a neighborhood of the organizing center at $(p, r) = (1, 1)$.
In this section, we present the details of the Bogdanov--Takens
unfolding procedure used to establish these results.
The section is somewhat technical, but since it is self-contained,
it can be skipped at first reading.

The unfolding is achieved by introducing a small positive
parameter~$\eta$ and rescaling the dependent and independent variables,
\begin{equation}
x (t) = \eta u (\tilde{t}) , \quad y (t) = \eta^2 v (\tilde{t}) , \quad \tilde{t} = \eta t .
\label{BT-2Dsym-unfolding}
\end{equation}
If $(x, y)$ is a solution of the system~(\ref{2Dsym-xy'}), 
then $(u, v)$ must satisfy the system
\begin{equation}
\begin{split}
\dot{u} &= v , \\
\dot{v} &= \mu u - u^3 + \eta (\lambda - u^2) v .
\end{split}
\label{2Dsym-uv,mu}
\end{equation}
Here, the dot $\dot{\ }$ stands for differentiation with respect to the variable~$\tilde{t}$,
and $\lambda$ and $\mu$ are parameters, which are defined in terms of $p$ and $r$,
\begin{equation}
\lambda = \frac{r - 1}{\eta^2}  , \quad \mu = \frac{r - p}{\eta^2} .
\label{2Dsym-lambda,mu}
\end{equation}
Note that $\mu$ is negative in region~I and positive in regions~II and~III 
(Figure~\ref{f-2Dsym-stability}). 
Henceforth, we omit the tilde and write $t$, instead of $\tilde{t}$.
\medskip

\noindent\textbf{Remark.}
The definition~(\ref{2Dsym-lambda,mu}) of $\lambda$ and $\mu$
generates a linear relation between~$p$ and~$r$,
\begin{equation}
(\lambda - \mu) (r-1) = \lambda (p-1) .
\label{2Dsym-pr-pencil}
\end{equation}
This is the equation of a pencil through the organizing center $(1,1)$
parameterized by $\lambda$.
Referring to the regions labeled I, II, and III in Figure~\ref{f-2Dsym-stability},
we note that $\lambda$ increases from 0 to infinity as one rotates counterclockwise
from the horizontal line $\{ p>1, r=1\}$ through region~I, 
then decreases as one continues to rotate through region~II, 
until $\lambda = 1$ at the vertical line $\{ p=1, r > 1 \}$,
and decreases further as one rotates through region~III, 
until $\lambda = 0$ at the horizontal line $\{ 0 < p < 1, r=1\}$.
\medskip

The results of the local analysis of Sections~\ref{ss-2Dsym-EquilStab} and~\ref{ss-2DsymHopf}
may be recovered directly from the system~(\ref{2Dsym-uv,mu}), as follows.
The origin $(0, 0)$ is an equilibrium state of~(\ref{2Dsym-uv,mu})
for all $\lambda$ and $\mu$, and if $\mu > 0$, there are two additional
equilibrium states, $(\pm \surd \mu, 0)$.
A linearization of~(\ref{2Dsym-uv,mu}) with $\mu < 0$ about $(0,0)$ 
shows that the real parts of the two eigenvalues pass through zero
at $\lambda = 0$, which corresponds to the line of
supercritical Hopf bifurcations $\{ p>1, r=1\}$.
Similarly, a linearization of~(\ref{2Dsym-uv,mu}) with $\mu>0$
about $(\pm \surd \mu,0)$ shows that the real parts of the two eigenvalues
pass through zero at $\lambda=\mu$, which corresponds to the line of 
subcritical Hopf bifurcations $\{p=1, r>1\}$.

\subsubsection{Hamiltonian Structures\label{sss-HamiltonianStructures}}
The equilibrium states of the system~(\ref{2Dsym-uv,mu})
are independent of $\eta$, so they persist as $\eta \to 0$.
In the limit, (\ref{2Dsym-uv,mu}) reduces to the Hamiltonian system
\begin{equation}
\begin{split}
\dot{u} &= v , \\
\dot{v} &= \mu u - u^3 .
\end{split}
\label{2Dsym-uv,mu,Hamiltonian}
\end{equation}
Closed orbits of this system are level curves of the Hamiltonian,
$H (u, v) = \tfrac12 v^2 - \tfrac12 \mu u^2 + \tfrac14 u^4$.
If $\mu < 0$, $H(u, v)$ reaches its minimum value 0 at the origin,
so the closed orbits are nested and surround the origin.
The more interesting case is $\mu > 0$,
where $H(u, v)$ reaches its minimum value at $(\pm \surd \mu, 0)$,
and $H(u, v) = 0$ at the saddle point at the origin.
We will analyze the case $\mu>0$ in detail and return to
the case $\mu< 0$ in Section~\ref{sss-LimitCycles,I}.
Recall that $\mu>0$ implies that $r>p$,
so the following analysis applies to the regions~II and~III
in Figure~\ref{f-2Dsym-stability}.
 
\begin{figure}[ht]
\begin{center}
\resizebox{0.6\textwidth}{!}{\includegraphics{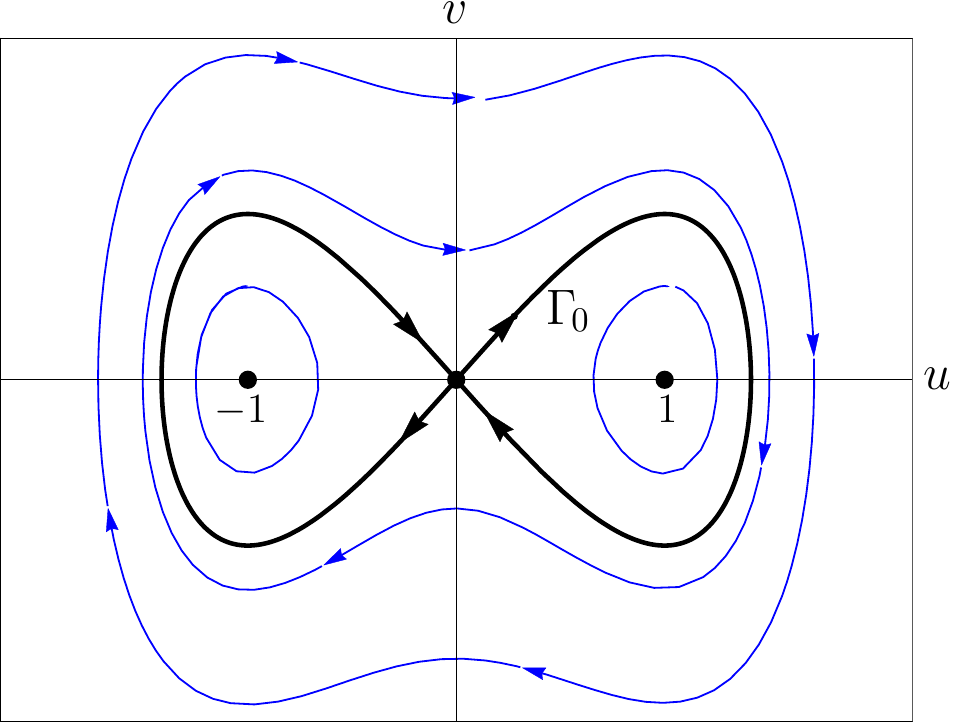}}
\caption{Phase portrait of the Hamiltonian system~(\ref{2Dsym-uv,mu,Hamiltonian}) with $\mu=1$.
 \label{f-2Dsym-PhasePortrait0}}
\end{center}
\end{figure}

Figure~\ref{f-2Dsym-PhasePortrait0} shows the phase portrait 
of the Hamiltonian system~(\ref{2Dsym-uv,mu,Hamiltonian}) with $\mu=1$.
(The phase portrait for other positive values of $\mu$ is similar.)
We see that there are several types of closed orbits.
There is a pair of homoclinic orbits to the origin,
there are periodic orbits in the interior of each of the homoclinic orbits
surrounding the equilibrium states $(\pm 1, 0)$, and 
there are large-amplitude periodic orbits 
external to the homoclinic orbits.
The question is, which of these closed orbits persist as the Hamiltonian
system~(\ref{2Dsym-uv,mu,Hamiltonian}) is perturbed to the system~(\ref{2Dsym-uv,mu}).
Because of the $\Z_2$ symmetry, it suffices to consider
the homoclinic orbit in the right half of the $(u, v)$ plane
and the periodic orbits in its interior;
the results for the closed orbits in the left half of the $(u,v)$ plane
follow by reflection.
Of course, we also need to consider the large-amplitude periodic orbits
that are external to the homoclinic orbits.
Notice the clockwise orientation of all these orbits.

\subsubsection{Melnikov Function\label{sss-Melnikov}}
The persistence of closed (homoclinic or periodic) orbits of Hamiltonian systems
under perturbations can be analyzed by means of the Melnikov function.
This function dates back at least to Poincar\'{e}~\cite{Poincare1890}; 
it features in articles by Melnikov~\cite{Melnikov1963} and Arnold~\cite{Arnold1964}
and in the book by Andronov et al.~\cite{Andronov1971}.
A definitive discussion can be found in the book of
Guckenheimer and Holmes~\cite[\S 4.5]{Guckenheimer1983}.
The Melnikov function and the associated theory apply to
a range of different systems, but for our purposes it suffices
to summarize the results for the general system
\begin{equation}
\begin{split}
\dot{u} &= v , \\
\dot{v} &= f(u) + \eta g(u, v) .
\end{split}
\label{BT-fg}
\end{equation}
(The functions $f$ and $g$ are not to be confused with
those in Section~\ref{ss-2DsymHopf}.)
In the limit as $\eta \to 0$, 
this system reduces to the Hamiltonian system
\begin{equation}
\begin{split}
\dot{u} &= v , \\
\dot{v} &= f(u)  .
\end{split}
\label{BT-Hamiltonian}
\end{equation}
Let $\Gamma_0 = \{ t \mapsto (u (t), v (t)), \, t \in I \}$ be any closed orbit
of~(\ref{BT-Hamiltonian}).
The \emph{Melnikov function} associated with $\Gamma_0$
is the integral $\int_I g (u(t), v(t)) v(t) \, dt$.
Thus, the Melnikov function measures the cumulative effect of the projection 
of the perturbed component of the vector field, $[ 0 \; g ]^\mathrm{t}$,
on the normal vector, $[ -f \; v ]^\mathrm{t}$, of the unperturbed vector field
along $\Gamma_0$.
If the Melnikov function vanishes on $\Gamma_0$,
then there exists---under suitable nondegeneracy conditions---a
family of closed orbits $\Gamma_\eta$ of the perturbed system~(\ref{BT-fg}) 
which are  $\mathcal{O}(\eta)$-close to $\Gamma_0$ as $\eta \to 0$.
Moreover, if the Melnikov function vanishes on $\Gamma_0$ and is positive (negative)
on nearby orbits that are to the right as $\Gamma_0$ is traversed, 
then the $\Gamma_\eta$ are locally stable (unstable).

\subsubsection{Dynamics in Regions~II and~III\label{sss-Dynamics,II,III}}
We apply the general results of the previous section
to the closed orbits of~(\ref{2Dsym-uv,mu}) identified 
at the end of Section~\ref{sss-HamiltonianStructures}:
either the homoclinic orbit in the right half of the $(u,v)$ plane,
or one of the periodic orbits in its interior, or one of the
large-amplitude periodic orbits external to the homoclinic orbits
(Figure~\ref{f-2Dsym-PhasePortrait0}).
We assume, without loss of generality, that $\mu = 1$,
so $f (u) = u - u^3$ and $g (u, v) = (\lambda - u^2) v$.

Let $\gamma = \{ (u(t), v(t)) : t \in \R \}$ be any of these closed orbits.
To indicate a particular orbit, we label $\gamma$ by the maximum value
of its first coordinate, $u (t)$, on its trajectory.
We consider this label as a variable and denote it by~$x$
(not to be confused with the dependent variable~$x$ in the Maasch-Saltzman model).
Thus, the function $\gamma : x \mapsto \gamma (x)$ is defined for all $x \in (1, \infty)$.
Specifically, $\gamma (x)$ is the homoclinic orbit in the right half plane if $x = \surd 2$,
a periodic orbit inside this homoclinic if $x \in (1, \surd 2)$, and
a large-amplitude periodic orbit that is external to the double homoclinic if $x > \surd 2$.

\textbf{Remark.} There are several ways to choose an identifier for $\gamma$.
For example, we could equally well have chosen the level-set value $h = H(\gamma)$,
as was done in~\cite{Kin2001}.
\medskip

Consider the closed orbit $\gamma (x) = \{ (u_x (t), v_x (t)) : t \in I(x) \}$ for any $x > 1$,
where $I (x) = \R$ if $x = \surd 2$, and $I(x)$ is a period interval otherwise. 
The Melnikov function associated with $\gamma (x)$ is
\be
M (\lambda, x) = \int_{I(x)} (\lambda - u^2(t))v^2(t) \, dt 
= \oint_{\gamma (x)} (\lambda - u^2 ) \, v (u) \, du = \lambda I_0 (x) - I_2 (x) ,
\label{2Dsym-Melnikov,x}
\ee
where $I_0$ and $I_2$ are defined by
\be
I_0 (x) = \oint_{\gamma (x)} v (u) \, du , \quad I_2 (x) = \oint_{\gamma (x)} u^2 \, v (u) \, du .
\label{2Dsym-I0I2}
\ee
Here, we have used the relation $v = \dot{u}$ to convert the time integral 
to a contour integral on the closed orbit $\gamma (x)$.

Recall that $\gamma (x)$ is oriented clockwise;
hence, Green's theorem yields the identity
$I_0 (x) = \iint_{D(x)} \, du \, dv$,
where $D(x)$ is the domain enclosed by $\gamma (x)$.
Consequently, $I_0 (x) > 0$, so the condition
$M (\lambda, x) = 0$ is satisfied if and only if
\be
\lambda = R (x) , \quad \mbox{where }R (x) = \frac{I_2 (x)}{I_0(x)} .
\label{2Dsym-R,a}
\ee
If, given $\lambda$, $M(\lambda, x) = 0$ at $x = \tilde{x}$,
then $M( \lambda, x) > 0$ for nearby orbits that are to the right
of $\gamma (\tilde{x})$ if $R' (\tilde{x}) < 0$,
and $M( \lambda, x) < 0$ for such nearby orbits if $R' (\tilde{x}) > 0$.
Hence, the local stability of closed orbits is determined by the sign of $R'(\tilde{x})$.

Since the components of the homoclinic and periodic orbits are known
in terms of hyperbolic and elliptic functions, respectively, 
$R(x)$ can be evaluated explicitly.
The computations were first done by Carr~\cite{Carr1981}
and subsequently refined by Cushman and Sanders~\cite{Cushman1985}.
For example, for the homoclinic orbit,
\be
\gamma (\surd 2) = \{ (\surd{2} \, \sech \, t, - \surd{2} \, \sech \, t \tanh \, t) , \; t \in \R \} ,
\ee
and $M(\lambda, \surd 2) = \tfrac43 \lambda - \tfrac{16}{15}$.
It follows that $M (\lambda, \surd 2) = 0$ if and only if $\lambda = \tfrac45$.
Moreover, $\lambda=\tfrac45$ is a simple zero.
Therefore, for all sufficiently small $\eta$ there exists
a $\lambda (\eta) = \tfrac45  + \mathcal{O} (\eta)$
and a homoclinic orbit near $\gamma (\surd 2)$, 
with a symmetric result in the left half of the $(u,v)$ plane.
In the $(p,r)$ plane, the homoclinic bifurcation curve is, to leading order,
tangent to the line $r-1 = -4 (p-1)$ at the organizing center $(1, 1)$; 
see~(\ref{2Dsym-pr-pencil}).

The figure below shows the graph of $R: x \mapsto R(x)$ for $1 < x < 2$. \\
\begin{minipage}{0.4\textwidth}
\resizebox{\textwidth}{!}{\includegraphics{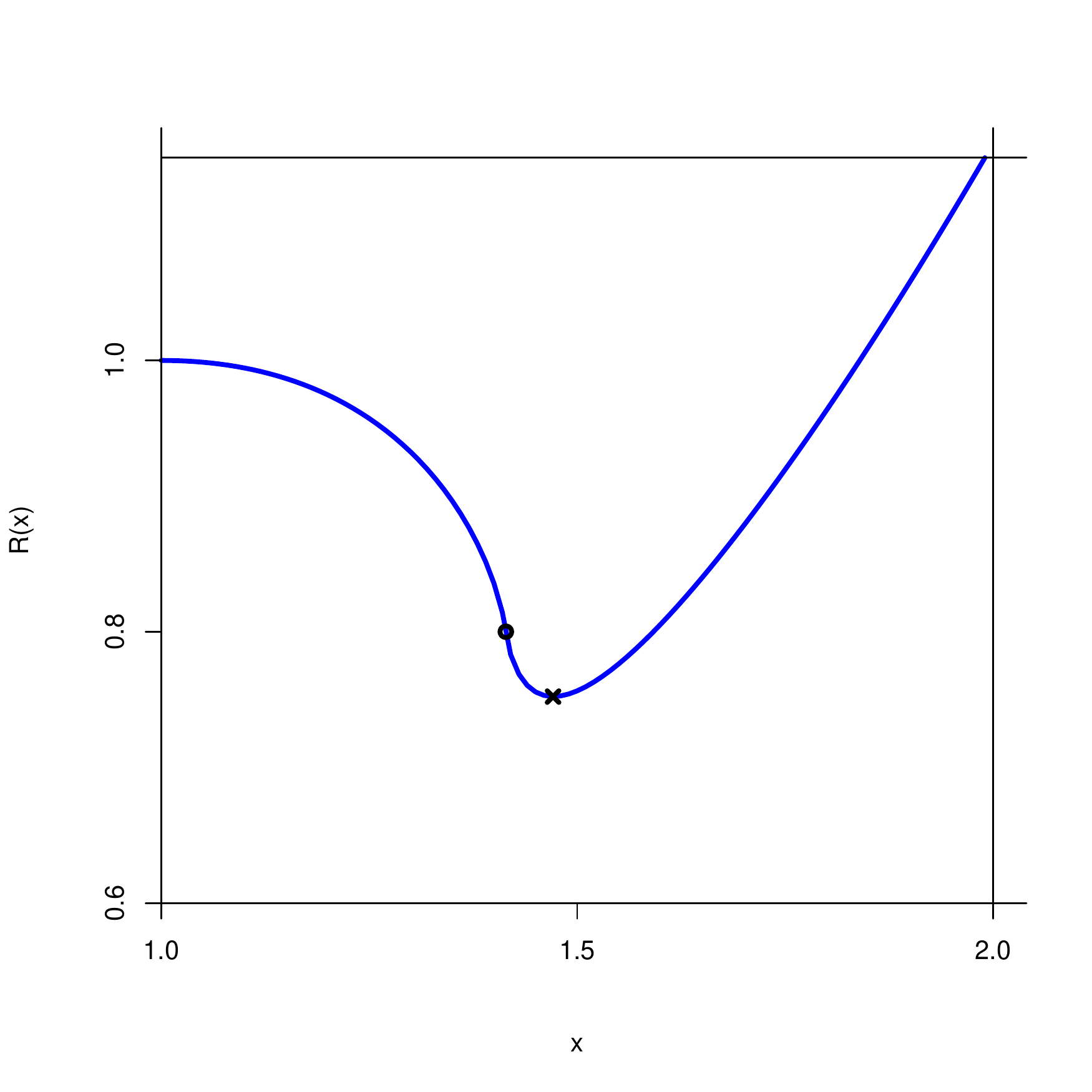}}
\end{minipage}
\hfill
\begin{minipage}{0.45\textwidth}
Starting from the values $R(1) = 1$ and $R'(1) = 0$,
$R (x)$ and $R' (x)$ decrease monotonically as $x$ increases.
At the homoclinic orbit (marked by a black dot),
$x = \surd 2$, $R(\surd 2) = \tfrac45$ and $\lim_{x \to \surd 2} R' (x) = - \infty$.
Beyond the homoclinic orbit, $R(x)$ decreases further, 
while $R' (x)$ increases until $R' (x) = 0$; 
at that point (marked by a black cross), 
$x = x^* \approx 1.471$ and $R'' (x^*) > 0$,
so $R(x)$ reaches its minimum value 
$R_{\min} (x^*) = \lambda^* \approx 0.752$.
Beyond this point, $R(x)$ increases monotonically; 
$R (x) \sim x^2$ as $x \to \infty$.
\end{minipage}

\vspace{0.4ex} 
\noindent  
Proofs of these statements, which do not use elliptic functions,
can be found, for example, in~\cite{Kin2001}.
It follows that closed orbits exist only for $\lambda > \lambda^*$,
and they are locally stable only for $x > x^*$.

\subsubsection{Limit Cycles in Regions~II and~III\label{sss-LimitCycles}}
The properties of the Melnikov function listed in the previous section
lead to the following results for the dynamics of the perturbed system~(\ref{2Dsym-uv,mu})
with~$\mu=1$.
In all statements, it is assumed that $\eta$ is sufficiently small positive.
\begin{itemize}
\item 
For $\lambda > 1$ (region~II), there are only stable large-amplitude limit cycles which
encircle the origin and pass through points $(x, 0)$ with $R(x) > 1$;
\item 
For $\frac{4}{5} < \lambda < 1$ (region~IIIa), there are stable large-amplitude limit cycles
which encircle the origin and unstable limit cycles in their interior which encircle
the equilibrium points $(\pm 1, 0)$;
\item 
For $\lambda = \frac{4}{5} + \mathcal{O}(\eta)$, there is a symmetric pair of unstable homoclinic orbits,
one in each half plane;
\item 
For $0.752\ldots < \lambda < \tfrac{4}{5}$ (region~IIIb), there is a stable large-amplitude limit cycle 
and an unstable large-amplitude limit cycle in its interior, both encircling the origin;
\item 
At $\lambda = 0.752\ldots + \mathcal{O}(\eta)$, the stable and unstable large-amplitude limit cycles
join in a saddle-node bifurcation;
\item 
For $\lambda < 0.752\ldots$ (region~IIIc), there are no limit cycles.
\end{itemize}

Thus, in addition to the homoclinic bifurcation curve found earlier,
which is tangent to the line $r-1 = -4 (p - 1)$ at the organizing center,
there is a curve of saddle-node bifurcations of limit cycles,
which, to leading order, is tangent to the line $r-1 \approx - 3.03 (p-1)$ at the organizing center.
This follows from~(\ref{2Dsym-pr-pencil}), with $\mu=1$ and $\lambda^* \approx 0.752$.
The two bifurcation curves partition the region~III of Figure~\ref{f-2Dsym-stability}
into the three regions~IIIa, IIIb, and IIIc, as sketched in Figure~\ref{f-2Dsym-stability2}
and superimposed on the color map of Figure~\ref{f-2Dsym-HeatPlot}.

\subsubsection{Limit Cycles in Region~I\label{sss-LimitCycles,I}}
It remains to investigate the dynamics of the system~(\ref{2Dsym-uv,mu})
for $\mu<0$ ($r < p$, region~I in Figure~\ref{f-2Dsym-stability}).
This case is considerably simpler than the case $\mu > 0$.
Without loss of generality, we may assume that $\mu = -1$.
The Hamiltonian is $H(u, v) = \tfrac12 v^2 + \tfrac12 u^2 + \tfrac14 u^4$.
The closed orbits can again be identified by the maximum value, $x$,
of its first coordinate~$u(t)$, which in this case ranges over all positive values, $x > 0$.
The Melnikov function is given by the same expression~(\ref{2Dsym-Melnikov,x})
and vanishes if $\lambda = R(x)$, as in~(\ref{2Dsym-R,a}).
In this case, both $R(x)$ and $R'(x)$ increase as $x$ increases,
so the Melnikov theory establishes that, for each $x>0$, 
there exists a value of $\lambda$ (given by the simple zero 
of the Melnikov function) such that,
for some $\lambda (\eta)$ that is $\mathcal{O}(\eta)$ close to this value,
the perturbed system has a unique limit cycle.

\section{The Asymmetric Two-Dimensional Model\label{s-2Dasym}}
Having a complete understanding of the dynamics of
the symmetric two-dimensional model~(\ref{2Dsym-xy}),
we are in a position to study the effects of symmetry breaking.
The asymmetric two-dimensional model is derived formally from
the Maasch--Saltzman model~(\ref{MS-xyz}) by setting
$q = \infty$ ($z = - x$),
\begin{equation}
\begin{split}
\dot{x} &= - x - y , \\
\dot{y} &= ry + px + s x^2 - x^2 y .
\end{split}
\label{2Dasym-xy}
\end{equation}
We will see that this system has two nondegenerate Bogdanov--Takens points,
which act as organizing centers in the $(p,r)$ parameter space.
The geometry of these organizing centers and the bifurcation curves
emanating from them may be understood naturally as a result of
the breaking of the lone $\Z_2$-symmetric Bogdanov-Takens point
studied in Section~\ref{ss-2DsymBT}.

\subsection{Equilibrium States and Their Stability\label{ss-2DasymEquilStab}}
\begin{figure}[ht]
\begin{center}
\resizebox{0.7\textwidth}{!}{\includegraphics{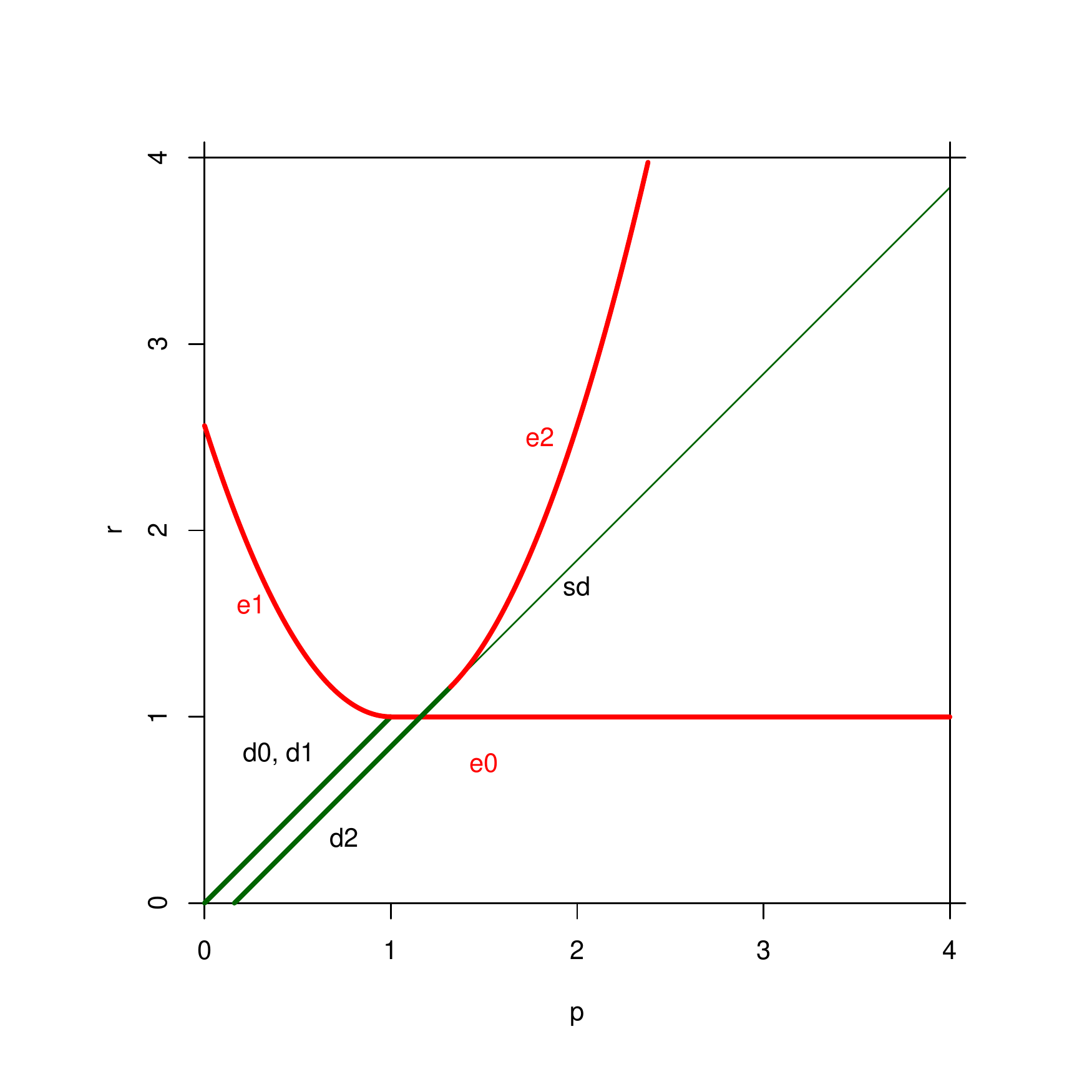}}
\caption{Stability regions of $P_0$, $P_1$, and $P_2$ for (\ref{2Dasym-xy}) with $s = 0.8$.
\label{f-2Dasym-stability}}
\end{center}
\end{figure}

The origin $P_0 = (0,0)$ is an equilibrium state of~(\ref{2Dasym-xy})
for all (positive) values of $p$, $r$, and $s$.
If $r > p - \tfrac14 s^2$, there are two additional equilibrium states,
namely $P_1 =  (x_1^*, - x_1^*)$ and $P_2 =  (x_2^*, - x_2^*)$,
where
\begin{equation}
x_1^* = \tfrac12 [- s + \sqrt{s^2 + 4(r - p)}] , \quad
x_2^* = \tfrac12 [- s - \sqrt{s^2 + 4(r - p)}] .
\label{2Dasym-x1*x2*}
\end{equation}
We refer to the line $r = p - \tfrac14 s^2$ as the \emph{shifted diagonal}
(marked `sd' in Figure~\ref{f-2Dasym-stability}). 
Note that $x_2^* < x_1^* < 0$ if $p - \tfrac14 s^2 < r < p$,
and $x_2^* < 0 < x_1^*$ if $r > p$.

Let $P = (x^*, -x^*)$ be any of the equilibrium states, with $x^* = 0$, $x_1^*$, or $x_2^*$.
The Jacobian of the vector field at $P$ is
\begin{equation}
\begin{pmatrix} -1 & -1 \\ p + 2s x^* + 2 (x^*)^2 & r - (x^*)^2 \end{pmatrix} .
\label{2Dasym-Jacobian}
\end{equation}
The system is linearly stable at $P$ if the trace is negative,
$-1 + r - (x^*)^2  < 0$,
and the determinant is positive,
$p - r + 3(x^*)^2 + 2s x^* > 0$.

The stability results are illustrated in Figure~\ref{f-2Dasym-stability}.
We see that the parameter space~$\R_+^2$ is partitioned
into six regions, which depend on $s$.
Referring to the labels in Figure~\ref{f-2Dasym-stability},
these regions are
\begin{equation}
\begin{split}
\mathrm{Oa} &= \left\{ (p, r) \in \R_+^2 : \mbox{ between `d0' and `d2', below `e0'} \right\} , \\
\mathrm{Ob} &= \left\{ (p, r) \in \R_+^2 : \mbox{ right of `d2', below `e0'} \right\} , \\
\mathrm{I} &= \left\{ (p, r) \in \R_+^2 : \mbox{ right of `sd', above `e0'}  \right\} , \\
\mathrm{IIa} &= \left\{ (p, r) \in \R_+^2 : \mbox{ between `e2' and `sd'} \right\} , \\
\mathrm{III} &= \left\{ (p, r) \in \R_+^2 : \mbox{ left of `d0' and `e1'} \right\} , \\
\mathrm{IIIo} &= \left\{ (p, r) \in \R_+^2 : \mbox{ between `e1' and `e2', above `e0' and `d2'} \right\} .
\end{split}
\label{2Dasym-O-III}
\end{equation}
Summarizing the results of the stability analysis, we find that

\begin{itemize}
\item $P_0$ is stable in regions~Oa and~Ob,
undergoes a supercritical Hopf bifurcation on~`e0' 
with natural frequency $\omega_0^* = \sqrt{p-1}$;
\item $P_1$ is stable in region~III,
undergoes a subcritical Hopf bifurcation on~`e1'
with natural frequency $\omega_1^* = \sqrt{2(r-1) + s \sqrt{r-1}}$; and
\item $P_2$ is stable in regions~Oa, III, and~IIIo,
undergoes a subcritical Hopf bifurcation on~`e2'
with natural frequency $\omega_2^* =  \sqrt{2(r-1) - s \sqrt{r-1}}$.
\end{itemize}

The introduction of asymmetry ($s > 0$) results in two changes.
The vertical line $\{ p=1, r>1 \}$, 
which is the locus of Hopf bifurcations for $P_1$ and $P_2$
in the symmetric case (Figure~\ref{f-2Dsym-stability}),
unfolds into a parabola.
The vertex of this parabola is at the point $(1, 1)$,
and the parabola is tangent to the shifted diagonal $r = p - \tfrac14 s^2$
at the point $(1 + \tfrac12 s^2, 1 + \tfrac14 s^2)$.
As we will see in Section~\ref{ss-2DasymBT},
both these points are organizing centers.
For convenience, we label them
\begin{equation}
Q_1 = (1, 1) , \quad 
Q_2 = (1 + \tfrac12 s^2, 1 + \tfrac14 s^2) .
\label{Q1Q2}
\end{equation}
As $(p, r)$ moves across the shifted diagonal in the direction
of decreasing $p$, the equilibrium states $P_1$ and $P_2$
emerge in a saddle-node bifurcation.
If the crossing occurs above $Q_2$, $P_1$ and $P_2$ are both unstable;
if it occurs below $Q_2$, $P_1$ is unstable while $P_2$ is stable.
In the region~Oa, $P_0$ and $P_2$ co-exist as stable equilibria,
while $P_1$ is an unstable equilibrium.
On the diagonal for $p < 1$, $P_0$ and $P_1$
exchange stability in a transcritical bifurcation.

\subsection{Organizing Centers\label{ss-2DasymBT}}
We make the change of variables $(x, -(x+y)) \mapsto (x, y)$
as in Section~\ref{ss-2DsymBT}.
In the new coordinate system, (\ref{2Dasym-xy}) becomes
\begin{equation}
\begin{split}
\dot{x} &= y , \\
\dot{y} &= (r - p) x + (r - 1) y - (s + y) x^2 - x^3 .
\end{split}
\label{2Dasym-xy'}
\end{equation}
The equilibrium points are
$P_0 = (0, 0)$, $P_1 = (x_1^*, 0)$, and $P_2 = (x_2^*, 0)$,
where $x_1^*$ and $x_2^*$ are again given by~(\ref{2Dasym-x1*x2*}).

Let $P = (x^*, 0)$ be any of the equilibrium points.
The Jacobian of the vector field at $P$ is
\begin{equation}
\begin{pmatrix} 0 & 1§ \\ r - p - 2s x^* - 3 (x^*)^2  & r - 1 - (x^*)^2 \end{pmatrix} .
\label{2Dasym-Jacobian'}
\end{equation}
If $x^* = 0$, the Jacobian has a double-zero eigenvalue at $Q_1$ for any~$s$,
and if $x^* = x_1^*$ or $x^* = x_2^*$, it has a double-zero eigenvalue at $Q_2$.
Hence, the introduction of asymmetry causes the organizing center
to unfold into a center at $Q_1$ associated with $P_0$
and a center at $Q _2$ associated with $P_1$ and $P_2$.

\begin{figure}[ht]
\begin{center}
\resizebox{0.7\textwidth}{!}{\includegraphics{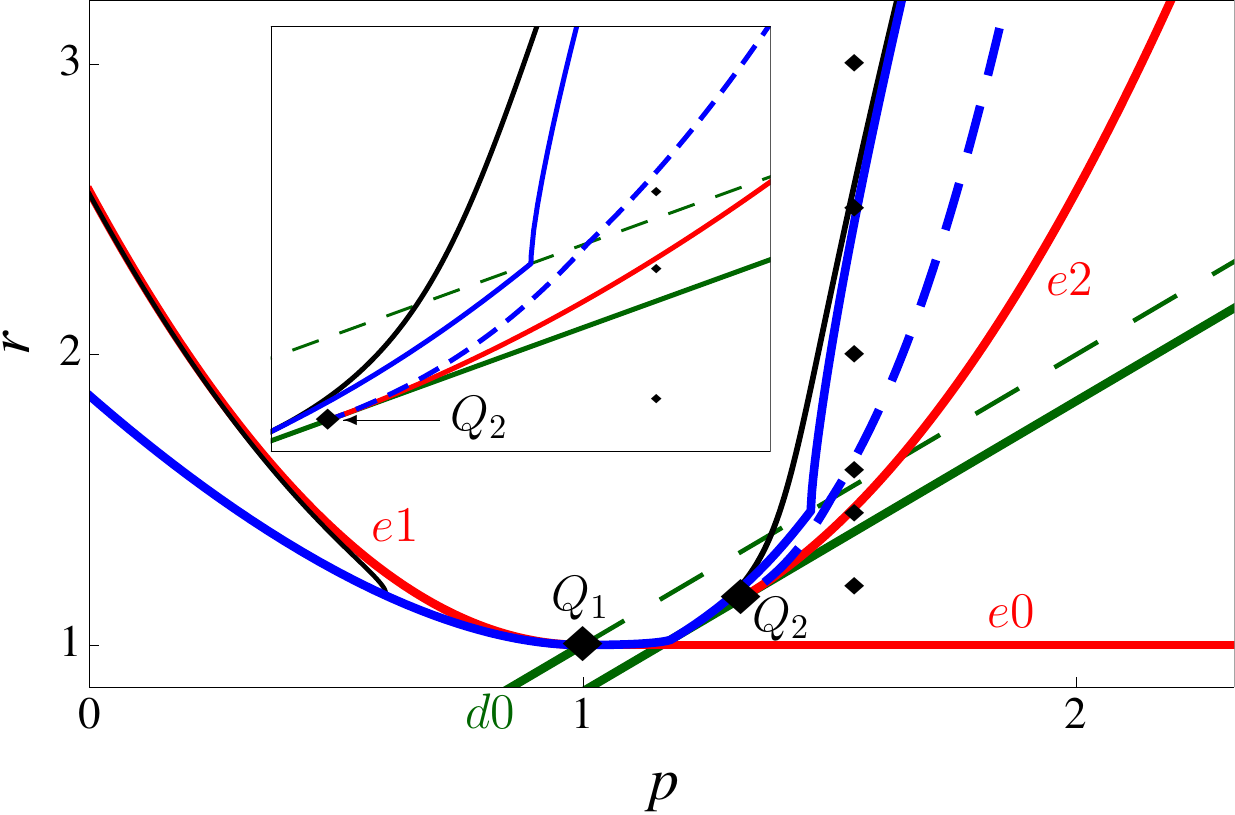}}
\caption{Stability boundaries and bifurcation curves for the system~(\ref{2Dasym-xy}) with $s = 0.8$.
The inset shows a neighborhood of $Q_2$.
\label{f-2Dasym-BifurcationCurves}}
\end{center}
\end{figure}

Figure~\ref{f-2Dasym-BifurcationCurves} shows the bifurcation curves
emanating from the two organizing centers for $s=0.8$
(the value chosen by Maasch and Saltzman).
They were computed with the AUTO continuation package.

There are three Hopf bifurcation curves (shown in red),
two emanating from $Q_1$ and one emanating from $Q_2$:
(i)~a parabolic curve to the left of $Q_1$, 
where $P_1$ undergoes a subcritical Hopf bifurcation,
(ii)~a horizontal line $\{r=1, p>1\}$ to the right of $Q_1$, 
where $P_0$ undergoes a supercritical Hopf bifurcation;
(iii)~a parabolic curve to the right of $Q_2$,
along which $P_2$ undergoes a subcritical Hopf bifurcation.
These three curves are the same as the ones identified in
the local analysis, cf.\ Section~\ref{ss-2DasymEquilStab}.
In addition, there are three homoclinic bifurcation curves (shown in blue),
two emanating from $Q_1$ (solid blue) and one emanating from $Q_2$ (dashed blue).
These curves are identified by an unfolding procedure,
as in Section~\ref{ss-BTUnfolding}.
Since the points $Q_1$ and $Q_2$ are both nondegenerate Bogdanov--Takens points,
there are no other bifurcation curves besides the Hopf and homoclinic bifurcation curves
emanating from them.

\begin{figure}[ht]
\begin{center}
\resizebox{0.9\textwidth}{!}{\includegraphics{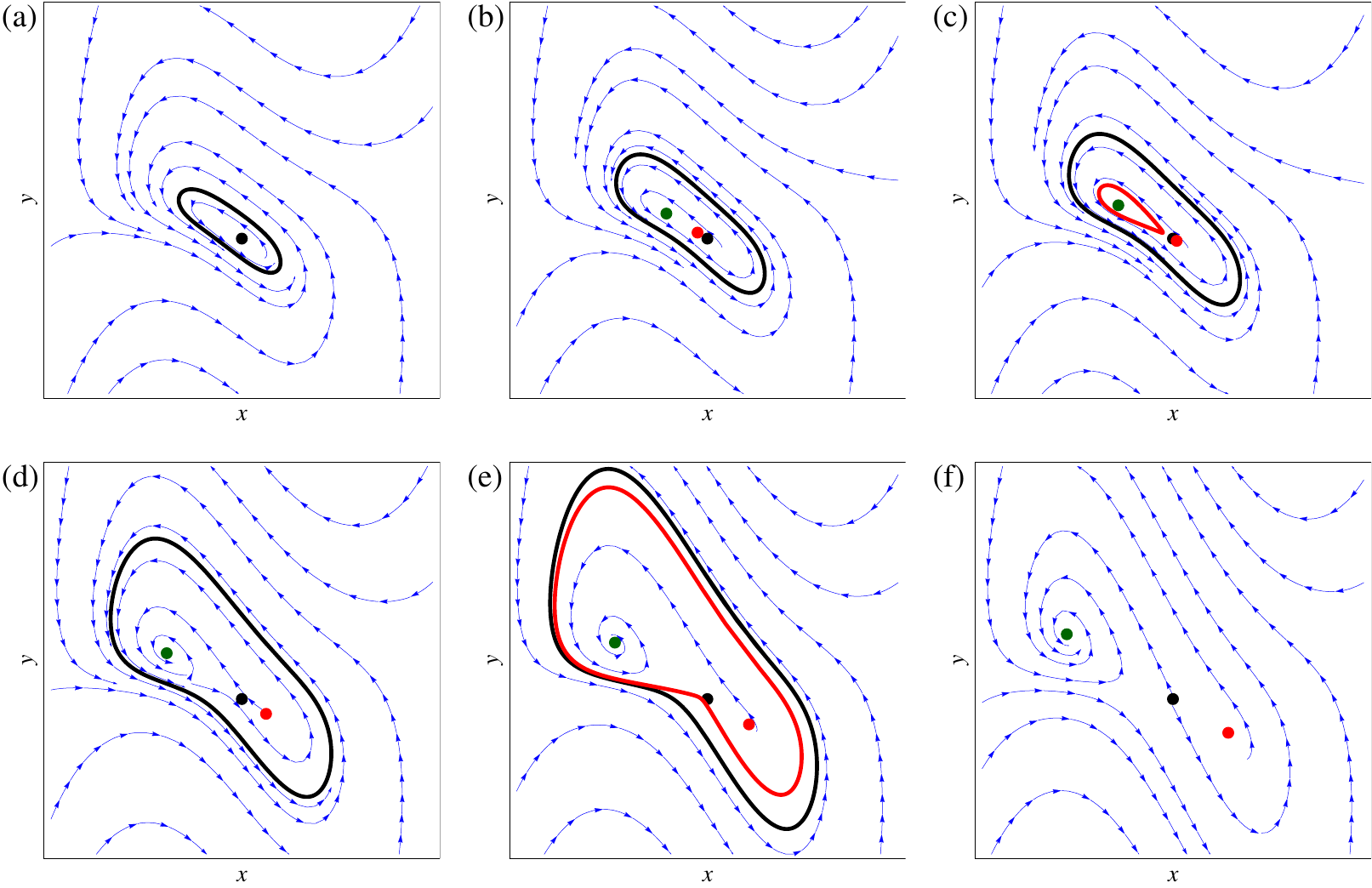}}
\caption{Phase planes of the system~(\ref{2Dasym-xy})
at the six small, black diamond markers in Figure~\ref{f-2Dasym-BifurcationCurves}.
The markers all lie on the vertical line $p=1.55$;
(a)~$r=1.2$, (b)~$r=1.45$, (c)~$r=1.6$, (d)~$r=2.0$, (e)~$r=2.5$, (f)~$r=3.0$.
\label{f-2Dasym-PhasePlanes}}
\end{center}
\end{figure}

Figure~\ref{f-2Dasym-PhasePlanes} shows the phase portraits at the six small black diamond
markers along the vertical line $p = 1.55$ in Figure~\ref{f-2Dasym-BifurcationCurves}.
The color scheme is as follows:
The flow of~(\ref{2Dasym-xy}) is shown as blue streamlines.
The black and red curves correspond to the stable and unstable limit cycles, respectively.
The equilibrium states $P_0$, $P_1$, and $P_2$ are indicated by
black, red, and green markers, respectively.
\begin{itemize}
\item Frame~(a), $r=1.2$: a stable limit cycle surrounds the unstable equilibrium point~$P_0$.
\item Frame~(b), $r=1.45$: the unstable equilibrium points $P_1$ and $P_2$ exist but are unstable, a stable limit cycle surrounds $P_1$ and $P_2$.
\item Frame~(c), $r=1.6$: the equilibrium points $P_1$ and $P_0$ have switched positions, $P_2$ has become stable, an unstable limit cycle surrounds $P_2$.
\item Frame~(d), $r=2.0$: the unstable limit cycle has disappeared in the homoclinic bifurcation.
\item Frame~(e), $r=2.5$: a large-amplitude unstable limit cycle exists inside the large-amplitude stable limit cycle.
\item Frame~(f), $r=3.0$: the large-amplitude stable and unstable limit cycles have disappeared in a saddle-node bifurcation; $P_2$ is the only attractor.
\end{itemize}

\section{Summary\label{s-Summary}}
The purpose of this chapter is to show an interesting application
of dynamical systems theory to a problem of climate science.
The object of investigation is a model of the Pleistocene climate,
proposed by Maasch and Saltzman in 1990.
The model  is a \emph{conceptual} model designed specifically
to explain the persistence of glacial cycles during the Pleistocene Epoch.
The Milankovitch theory of orbital forcing establishes a correlation
between glacial cycles and periodic oscillations in the Earth's orbit
around the Sun, but orbital forcing by itself is insufficient
to explain the observed temperature changes.

The Maasch--Saltzman model incorporates a feedback mechanism
that is driven by greenhouse gases, in particular atmospheric \CO.
The model is based on plausible physical arguments, and preliminary
computational experiments indicate that it reproduces several salient
features of the Pleistocene temperature record.
In this article, the emphasis is on the internal dynamics of the model.
We focused on the prevalence and bifurcation properties of limit cycles
in the various parameter regimes in the absence of external forcing.

The Maasch--Saltzman model is formulated in terms of the anomalies
of the total global ice mass, the atmospheric \CO\ concentration, and
the volume of the North Atlantic Deep Water (NADW, a measure of
the strength of the North Atlantic overturning circulation).
It consists of three differential equations with four parameters and is
difficult to analyze directly.
Our results indicate how one can obtain fundamental insight into its complex
dynamics by first considering a highly simplified two-dimensional version.
The dimension reduction is achieved by (formally) letting one of the 
parameters---representing the rate of change of the volume of NADW
relative to that of the total global ice mass---tend to infinity.
The approximation is justified by the observation that the NADW
changes on a much faster time scale than the total global ice mass.

The two-dimensional model has two primary parameters, $p$ and~$r$,
which are both positive, and one secondary parameter~$s$, 
which reflects the asymmetry between the glaciation and deglaciation
phases of the glacial cycles.
By first ignoring the asymmetry $(s=0)$, we obtained a complete
understanding of the dynamics and the persistence of limit cycles.

Figure~\ref{f-2Dsym-stability2}, which is a sketch of the various bifurcation curves
in the $(p, r)$ parameter space, summarizes the main results.
The origin is an equilibrium state $P_0$ for all $(p, r)$; $P_0$ is stable in region~O.
In addition, there are two equilibrium states~$P_1$ and $P_2$ in regions~II and~IIIa-c;
they are generated in a pitchfork bifurcation along the diagonal $r=p$
and are stable in region~IIIa-c.
Stable limit cycles exist in regions~O, I, II, IIIa, and~IIIb. 
They are created in supercritical Hopf bifurcations along the boundary
between regions~O and~I.
In region~I, they surround the unstable equilibrium state $P_0$,
and in the other regions they surround all three equilibrium states.
Along the boundary between regions~II and~IIIa, a pair of unstable
limit cycles, each surrounding one of the two stable equilibrium 
states~$P_1$ and~$P_2$, are created in a subcritical Hopf bifurcation.
These newborn limit cycles grow in amplitude until they become
homoclinic orbits at the boundary between regions~IIIa and~IIIb;
in region~IIIb, they have merged into one large-amplitude unstable limit cycle,
which surrounds the three equilibrium states and sits just inside
the large-amplitude stable limit cycle.
Finally, the stable and unstable large-amplitude limit cycles merge
and disappear in a saddle-node bifurcation as $(p, r)$ crosses
the lower boundary of region~IIIb into region~IIIc.
All these results have been confirmed computationally;
see Figure~\ref{f-2Dsym-HeatPlot}.

Having thus obtained a complete understanding of the dynamics 
of the symmetric simplified model, we then re-introduced asymmetry ($s > 0$).
A comparison of Figure~\ref{f-2Dasym-BifurcationCurves} with
Figure~\ref{f-2Dsym-stability2} shows the effects of symmetry breaking.
The main change is that the single organizing center, which governed
the dynamics in the symmetric case, splits into two organizing centers.
Also, the curves of homoclinic bifurcations and saddle-node bifurcations
of limit cycles become more complex.

The complexity of the bifurcation diagram of Figure~\ref{f-2Dasym-BifurcationCurves}
gives an indication why it is difficult to analyze the dynamics of the Maasch--Saltzman model
directly.
It also justifies our approach of first analyzing the highly simplified model
and then gradually relax the constraints that were imposed
to derive the simplified model~\cite{EnglerKKVo2017}.

\end{document}